\documentclass{amsart}

\theoremstyle{definition}

\theoremstyle{remark}

\numberwithin{equation}{section}

\newcommand{\abs}[1]{\lvert#1\rvert}

\reversemarginpar

\newcommand{\nn}{\nonumber}
\newcommand{\no}{\noindent}
\newcommand{\cl}{\mathop{\rm Cl}\nolimits}
\newcommand{\li}{\mathop{\rm Li}\nolimits}
\newcommand{\realpart}{\mathop{\rm Re}\nolimits}
\newcommand{\imagpart}{\mathop{\rm Im}\nolimits}
\newcommand{\lp}{\ln \sqrt{2 \pi}}
\newcommand{\ba}{\begin{eqnarray}}
\newcommand{\ea}{\end{eqnarray}}
\newcommand{\ift}{\int_{0}^{\infty}}
\newcommand{\ione}{\int_{0}^{1}}
\newcommand{\npR}{\mathbb{R}^{-}_{0}}
\newcommand{\nnR}{\mathbb{R}^{+}_{0}}
\newcommand{\allR}{\mathbb{R}}
\newcommand{\allN}{\mathbb{N}}
\newcommand{\nnN}{\mathbb{N}_{0}}
\newcommand{\st}{}
\newcommand{\stoli}{}
\newcommand{\stboth}{}

\begin{document}

\title[Hurwitz zeta function] {On some definite integrals involving the Hurwitz
zeta function}

\author{Olivier Espinosa}
\address{Departamento de F\'{\i}sica,
Universidad T\'{e}cnica Federico Santa Mar\'{\i}a, Valpara\'{\i}so, Chile}
\email{espinosa@fis.utfsm.cl}

\author{Victor H. Moll}
\address{Department of Mathematics,
Tulane University, New Orleans, LA 70118}
\email{vhm@math.tulane.edu}

\subjclass{Primary 33}

\date{\today}


\begin{abstract}
We establish a series of integral formulae involving the Hurwitz zeta
function. Applications are given to integrals of Bernoulli polynomials,
$\log \Gamma(q)$ and $\log \sin (q)$.
\end{abstract}

\maketitle


\newtheorem{Definition}{\bf Definition}[section]
\newtheorem{Thm}[Definition]{\bf Theorem}
\newtheorem{Lem}[Definition]{\bf Lemma}
\newtheorem{Cor}[Definition]{\bf Corollary}
\newtheorem{Prop}[Definition]{\bf Proposition}
\newtheorem{Example}[Definition]{\bf Example}

\section{Introduction} \label{S:intro}

The Hurwitz zeta function, defined by
\ba
\zeta(z,q) & = & \sum_{n=0}^{\infty} \frac{1}{(n+q)^{z}}
\stboth
\label{hurdef}
\ea
\no
for $z \in \mathbb{C}$ and $q \neq 0, \, -1, \, -2, \cdots,$ is one of the
fundamental transcendental functions. The series
converges for $\realpart(z) > 1$ so that $\zeta(z,q)$ is an
analytic function of $z$ in this region. The integral representation
\ba
\zeta(z,q) & = & \frac{1}{\Gamma(z)} \ift
\frac{e^{-qt}}{1 - e^{-t}} t^{z-1} dt,
\stboth
\label{intrep0}
\ea
\no
where $\Gamma(z)$ is Euler's gamma function, is valid
for $\realpart(z) > 1$ and $\realpart(q) > 0$, and can be used to show that
$\zeta(z,q)$ admits an analytic extension to the whole
complex plane except for a simple pole at $z=1$.
In most of the examples discussed here
we consider only the range $0 < q \leq 1$.
Special cases of $\zeta(z,q)$
include the Riemann zeta function
\ba
\zeta(z,1) & = & \zeta(z) = \sum_{n=1}^{\infty} \frac{1}{n^{z}}
\label{riezeta0}
\stboth
\ea
\no
and
\ba
\zeta(z, \tfrac{1}{2}) & = & 2^{z} \sum_{n=0}^{\infty} \frac{1}{(2n+1)^{z}}
= (2^{z}-1) \zeta(z).
\stboth
\label{riezeta1}
\ea

The function $\zeta(z,q)$ admits several integral representations in
addition to (\ref{intrep0}).
For example, Hermite proved
\ba
\zeta(z,q) & = &  \frac{1}{2}q^{-z} + \frac{1}{z-1} q^{1-z} +
2q^{1-z} \ift \frac{\sin( z \tan^{-1} t) dt}{(1+t^2)^{z/2} \left( e^{2 \pi t q}
-1 \right) },
\stboth
\label{repint}
\ea
\no
which is valid for $q > 0$ and $z \neq 1$. In fact, (\ref{repint}) is
an explicit
representation of the analytic continuation of (\ref{hurdef}) to
$\mathbb{C} - \{ 1 \}$.

Among the many places in which $\zeta(z,q)$ appears we mention the evaluation
by K\"{o}lbig \cite{kolbig1} of integrals of the form
\ba
R_{m}(\mu,\nu) & = & \ift e^{-\mu t} t^{\nu -1} \log^{m}t \, dt, \label{kol}
\ea
\no
an example of which is
\ba
 R_{2}(\mu,\nu) & = &
\mu^{-\nu} \Gamma(\nu) \left[ \left( \psi(\nu) - \ln \mu \right)^{2} +
\zeta(2, \nu) \right]. \label{kolbig11}
\ea
\no
Here
\ba
\psi(x)  & = & \Gamma'(x)/\Gamma(x) \label{digamma}
\ea
\no
is the logarithmic derivative of
$\Gamma(x)$, also called the digamma function. \\

The Hurwitz zeta function also plays a role in
Vardi's evaluations \cite{vardi2}
\ba
\int_{\pi/4}^{\pi/2} \log \log \tan x dx & = & \frac{\pi}{2}
\ln \left( \frac{\Gamma(3/4) \sqrt{2 \pi}}{\Gamma(1/4)} \right) \label{var}
\ea
\no
and \cite{vardi1} of Kinkelin's constant
\ba
\ln A & := & \lim\limits_{k \to \infty} \left[ \ln(1^{1} 2^{2} \cdots k^{k})
- \tfrac{1}{2} ( k^{2} + k + 1) \ln k + \tfrac{k^{2}}{4} \right]  \label{var2}
\ea
\no
as
\ba
\ln A & = & \text{exp}\left( \tfrac{1}{12} - \zeta'(-1) \right).  \label{var3}
\ea

Yue and Williams \cite{yuewil1, yuewil2} established the integral
representation
\ba
\zeta(z,q) & = & 2 (2 \pi)^{z-1} \ift \frac{ e^{x} \sin( \pi z/2 + 2 \pi q) -
\sin( \pi z/2)}{e^{2x} - 2e^{x} \cos 2 \pi q + 1} x^{-z} dx
\ea
\no
and used it to evaluate definite integrals, (\ref{var}) among
them.  For example, for $0 < a < 1$, they obtain
\ba
\int_{0}^{\infty} \frac{e^{-x} \ln x \; dx}{e^{-2x} -2e^{-x} \cos 2 \pi a + 1}
& = & \frac{\pi}{2 \sin 2 \pi a} \ln\left( \frac{\Gamma(1-a)}{(2 \pi)^{2a-1} \,
\Gamma(a) } \right). \nn
\ea

\medskip

Integrals involving the Hurwitz zeta function also appear in problems
dealing with distributions of
$\{ nx \}$ for $x \not \in
\mathbb{Q}$ and $n \in \allN$, where $\{ x \}$ denotes the fractional part
of $x$. In this context Mikolas \cite{mikolas}
established the identity
\ba
\ione \zeta(1-z,\{ a q \}) \, \zeta(1-z, \{ b q \}) dq & = &
2 \Gamma^{2}(z) \frac{\zeta(2z)}{( 2 \pi)^{2z}}
\left( \frac{(a,b)}{[a,b]} \right)^{z} \label{miko}
\ea
\no
for $a, \, b \in \allN$. Here $(a,b)$
is the greatest common divisor of $a$ and $b$ and $[a,b]$ is their least
common multiple. \\

\medskip

The Hurwitz zeta function also plays a role in the evaluation of
functional determinants that appear in mathematical physics. See
\cite{elizalde} for a miscellaneous list of physical examples.
The Hurwitz zeta function has also recently appeared in connection
with the problem of a gas of non-interacting electrons in the
background of a uniform magnetic field \cite{USM-TH-85}. For
instance, it is shown there that the {\em density of states}
$g(E)$, in terms of which all thermodynamic functions are to be
computed, can be written as
\ba\label{gE1}
g(E) = V\frac{{4\pi }}{{h^3 }}(2e\hbar B)^{1/2} E\,h_{1/2} \left(
{\frac{{E^2  - M^2 }}{{2e\hbar B}}} \right),
\ea
where $V$ stands for volume, $B$ for magnetic field, $M$ is the
electron mass, and
\ba\label{def-h}
h_{1/2} (q): = \zeta (\tfrac{1}{2},\left\{ q \right\}) - \zeta
(\tfrac{1}{2},q + 1) - \frac{1}{2q^{1/2}}.
\ea
As before, $\{ q \}$ in (\ref{def-h}) denotes the fractional part
of $q$.

\medskip

General information about $\zeta(z,q)$ appears in
\cite{bateman} and \cite{weis}.
\medskip

In this paper we derive a series of formulae for definite integrals containing
$\zeta(z,q)$ in the integrand. A search of the
standard tables of integrals reveals very few examples in \cite{pr} and none
in \cite{gr}. For
instance, in  \cite{pr}, section $1.2.1$  we find the
indefinite integral
\ba
\int \zeta(z,q) dq & = & \frac{1}{1-z} \zeta(z-1,q), \label{intind}
\stboth
\ea
\no
which is an elementary  consequence of
\ba
\frac{\partial}{\partial q} \zeta(z-1,q) & = & (1-z) \zeta(z,q).
\stboth
\label{partial}
\ea
\no
Section $2.3.1$ of \cite{pr} gives two definite integrals:
\ba
\ift q^{\alpha -1} \zeta(z,a + bq) dq & = & b^{-\alpha} B( \alpha, z - \alpha)
\zeta(z - \alpha, a)
\stboth
\nn
\ea
\no
for $a, b \in \mathbb{R}^{+}, \; 0 < \realpart(\alpha) < \realpart(z) -1;$ and
\ba
\ift q^{\alpha - 1} \left[ \zeta(z,q) - q^{-z} \right] dq & = &
B(\alpha, z - \alpha) \zeta(z- \alpha)
\stboth
\nn
\ea
\no
for $0 < \realpart(\alpha) < \realpart(z) - 1$, where $B(x,y)$ is
the beta function. The second integral is actually a special case
of the first with $a=b=1$.
The only other example in
\cite{pr} is the  evaluation of one of the Fourier coefficients of
$\zeta(z,q)$ in section $2.3.1$:
\ba
\ione \sin (2 \pi q) \zeta(z,q) dq & = & \frac{ (2 \pi)^{z}}{4 \Gamma(z)}
\, \text{ csc} \left(\frac{z \pi}{2} \right)
\stboth
\label{firstfou}
\ea
\no
for $1 < \realpart(z) < 2$.  \\

The tables \cite{gr, pr} do contain many examples involving the special case
\ba
\zeta(1-m,q) & = & - \frac{1}{m} B_{m}(q)
\stboth
\label{bernoulli}
\ea
\no
for $m \in \allN, \, q \in \nnR$, where $B_{m}(q)$ are the Bernoulli
polynomials defined by their generating function
\ba
\frac{t e^{qt}}{e^{t} -1} & = &  \sum_{m=0}^{\infty} B_{m}(q) \frac{t^{m}}{m!}
\stboth
\label{gener}
\ea
\no
for $|t| < 2 \pi$. These polynomials can be expressed as
\ba
B_{m}(q) & = & \sum_{k=0}^{m} \binom{m}{k} B_{k} q^{m-k} \label{polyber}
\stboth
\ea
\no
in terms of the Bernoulli numbers $B_{m} = B_{m}(0)$. The latter are rational
numbers; for example, $B_{0}=1, \; B_{1} = -1/2$,  and,
$B_{2} = 1/6$.  The Bernoulli numbers of odd index $B_{2m+1}$
 vanish for $m \geq 1$, and  those with even index
satisfy $(-1)^{m+1}B_{2m} > 0$. \\

The relation (\ref{polyber}) can be inverted to produce
\ba
q^{n} & = & \frac{1}{n+1}
\sum_{j=0}^{n} \binom{n+1}{j} B_{j}(q),
\stboth
 \label{polyinv}
\ea
\no
and since $B_{j}(1-q) = (-1)^{j}B_{j}(q)$, we also have

\ba
(1-q)^{n} & = & \frac{1}{n+1}
\sum_{j=0}^{n} (-1)^{j} \binom{n+1}{j} B_{j}(q).
\stboth
 \label{polyinv1}
\ea
\no
For example,
\ba
B_{0}(q) = 1, \;\;\; B_{1}(q) = q - \tfrac{1}{2}, \;\;\;
B_{2}(q) = q^{2} -q + \tfrac{1}{6} \nn
\stboth
\ea
\no
yield
\ba
1 = B_{0}(q), \;\;\; q = B_{1}(q) + \frac{1}{2}B_{0}(q), \;\;\;
 q^{2} = B_{2}(q) + B_{1}(q) +
\frac{1}{3} B_{0}(q). \nn
\stboth
\ea

The results presented here are consequences of the
Fourier expansion of $\zeta(z,q)$:
\begin{equation}
\zeta(z,q) = \frac{2 \Gamma(1-z)}{(2 \pi)^{1-z}}  \times
\left( \sin\left( \frac{ \pi z}{2} \right) \sum_{n=1}^{\infty}
\frac{\cos( 2 \pi q n)}{n^{1-z} }  +
\cos\left( \frac{ \pi z}{2} \right) \sum_{n=1}^{\infty}
\frac{\sin( 2 \pi q n)}{n^{1-z} }  \right).
\stboth
\label{fouzeta}
\end{equation}
\no
This expansion, valid for $\realpart(z) <0$ and $0 < q < 1$,  is due to
Hurwitz and is derived
in \cite{ww}, page 268.  A proof of (\ref{fouzeta}) based upon the
representation\footnote{$\lfloor{x\rfloor}$ is the floor of  $x$.}
\ba
\zeta(z,q) & = & z \int_{-q}^{\infty} \frac{\lfloor{x\rfloor} -x +
\tfrac{1}{2}} {(x+q)^{z+1}} dx
\st
\label{int3}
\ea
appears in \cite{hurber}. The result
\ba
\ione \zeta(z,q) dq & = & 0,
\stboth
\label{vanishing}
\ea
\no
valid for $\realpart(z)<0$,
follows directly from the representation (\ref{fouzeta}).
Although the Fourier expansion is derived strictly for $\realpart(z)<0$, it
also holds for the boundary value $z=0$. We shall thus simply take
$z \in  \npR$ in most of the formulae presented below. \\

Our goal is to employ the representation (\ref{fouzeta})
to evaluate definite integrals containing $\zeta(z,q)$ in the integrand.
These evaluations can be seen as examples of the {\em Hurwitz transform}
defined by
\ba
\mathfrak{H}(f) & :=  & \ione  f(q) \zeta(z,q) dq.
\stboth
\label{hurtrans1}
\ea
\no
Properties of $\mathfrak{H}$ and its uses will be discussed elsewhere.   The
relation (\ref{bernoulli}) between Bernoulli polynomials and the Hurwitz
zeta function yields, for each evaluation of the Hurwitz transform, an
explicit formula for an integral of the type
\ba
{\mathfrak{B}}_{m}(f) & :=  & \ione  f(q) B_{m}(q) dq,
\stboth
\label{bertrans1}
\ea
\no
and by (\ref{polyinv}) the evaluation of the moments of the function $f$
\ba
{\mathfrak{M}}_{n}(f) & :=  & \ione  q^{n} f(q) dq.
\stboth
\label{momtrans1}
\ea

We have attempted to evaluate symbolically, using Mathematica 4.0 and/or
Maple V, each of the examples presented here. The few cases in which this
attempt was successful are so indicated.  \\

The relations
\ba
\zeta(2n) & = & \frac{ (-1)^{n+1} (2 \pi)^{2n} B_{2n}  }{2
(2n)!},\quad n\in\nnN,
\stboth
\label{berzeta}  \\
\zeta(1-n) & = & \frac{(-1)^{n+1} B_{n}}{n}, \quad n \in\allN,
\stboth
\label{berzeta1} \\
\zeta'(-2n) & = & (-1)^{n} \frac{(2n)! \, \zeta(2n+1)}{2 \, (2 \pi)^{2n} },
\quad n \in\allN,
\stboth
\label{berzeta2}  \\
\zeta'(0)  & = & -\lp,
\stboth
\label{berzeta3}
\ea
\no
and Riemann's functional equation
\ba
\zeta(1-s) & = & \frac{\zeta(s) (2 \pi)^{1-s} }{2 \Gamma(1-s) \sin(\pi s/2) }
\stboth \label{Riemann1} \\
  & = &  2 \cos \left(\frac{\pi s}{2} \right)
\frac{\zeta(s) \Gamma(s) }{(2 \pi)^{s}} \stboth
\label{Riemann2}
\ea
\no
will be used to simplify the integrals discussed below.  The
form (\ref{Riemann2}) follows from (\ref{Riemann1}) by use of
the reflection formula
\ba
\Gamma(x) \Gamma(1-x ) & = & \frac{\pi}{\sin \pi x}
\stboth
\label{reflec}
\ea
\no
for the gamma function. The basic relation between the beta and gamma
functions,
\ba
B(x,y) & = & \frac{\Gamma(x) \, \Gamma(y) }{\Gamma(x+y)},
\stboth
\label{betagamma}
\ea
\no
will also be employed throughout.  \\

\section{The Fourier expansion of $\zeta(z,q)$.} \label{first}

In this section we employ the Fourier expansion (\ref{fouzeta}) for
$\zeta(z,q)$ to evaluate definite integrals of the form
\ba
\mathfrak{H}(f) & :=  & \ione  f(q) \zeta(z,q) dq. \label{hurtrans}
\stboth
\ea
\no
The expansion is valid for $z \leq 0$. Section \ref{extzneg} discusses the
extension of some of these evaluations to the case $z >0$.  \\

We first record the Fourier coefficients of $\zeta(z,q)$. These can be read
directly from (\ref{fouzeta}).

\begin{Prop}
The Fourier coefficients of $\zeta(z,q)$ are given by
\ba
\int_{0}^{1} \sin(2 k \pi q) \zeta(z,q) dq & = & \frac{(2 \pi)^{z} \, k^{z-1}}
{4 \Gamma(z)} \text{csc} \left( \frac{z \pi}{2} \right)
\stboth
 \label{fou1}
\ea
\no
and
\ba
\int_{0}^{1} \cos(2 k \pi q) \zeta(z,q) dq & = & \frac{(2 \pi)^{z} \, k^{z-1}}
{4 \Gamma(z)} \text{sec} \left( \frac{z \pi}{2} \right).
\stboth
\label{fou2}
\ea
\end{Prop}
\begin{proof}
The orthogonality of the trigonometric functions and (\ref{fouzeta}) yield
\ba
\ione \sin(2 k \pi q) \zeta(z,q) dq & = & \frac{\Gamma(1-z)}
{(2 \pi k)^{1-z}} \,
\cos \left( \frac{\pi z}{2} \right).
\stboth
\ea
\no
Now use the reflection formula (\ref{reflec}) to
obtain (\ref{fou1}). The calculation of (\ref{fou2}) is similar.
\end{proof}

The theorem below reduces the evaluation of an integral of the type
considered here to the evaluation of a Dirichlet series formed with
the Fourier coefficients of the integrand.  The remainder of the paper
are applications of this result.

\begin{Thm}
\label{main}
Let $f(w,q)$ be defined for $q \in [0,1]$ and a parameter $w$. Let
\ba
f(w,q) & = & a_{0}(w) +
\sum_{n=1}^{\infty} a_{n}(w) \cos(2 \pi q n) + b_{n}(w) \sin( 2
\pi q n)
\label{fouf}
\stboth
\ea
\no
be its Fourier expansion, so that
\ba
a_{n}(w) & = & 2 \ione f(w,q) \cos(2 \pi q n) dq,  \quad n \, \geq
0,
\stboth
\label{coeff1} \\
b_{n}(w) & = & 2 \ione f(w,q) \sin(2 \pi q n) dq, \quad n \geq 1.
\stboth
\label{coeff2}
\ea
\no
Then, for $z\in\npR$,
\begin{align}
\ione f(w,q) \zeta(z,q) dq &=
\frac{\Gamma(1-z) }{(2 \pi)^{1-z}}
\left( \sin \left( \frac{\pi z}{2} \right)
\sum_{n=1}^{\infty} \frac{a_{n}(w)}{n^{1-z}} +
\cos \left( \frac{\pi z}{2} \right)
\sum_{n=1}^{\infty} \frac{b_{n}(w)}{n^{1-z}}  \right)
\stboth
\label{mainint}\\
\intertext{and}
\ione f(w,q) \zeta(z,1-q) dq &=
\frac{\Gamma(1-z) }{(2 \pi)^{1-z}}
\left( \sin \left( \frac{\pi z}{2} \right)
\sum_{n=1}^{\infty} \frac{a_{n}(w)}{n^{1-z}} - \cos \left(
\frac{\pi z}{2} \right) \sum_{n=1}^{\infty}
\frac{b_{n}(w)}{n^{1-z}}  \right).
\label{mainint2}
\stboth
\end{align}
\end{Thm}
\begin{proof}
Multiply (\ref{fouf}) by $\zeta(z,q)$, integrate over $[0,1]$,  and apply
(\ref{fou1}) and (\ref{fou2}) to give (\ref{mainint}).  Observe that
the integral
of $\zeta(z,q)$ over $[0,1]$ vanishes, so there is no contribution from
$a_{0}(w)$. The second result follows from the fact that the
Fourier expansion of $\zeta(z,1-q)$ differs from that of
$\zeta(z,q)$ given in (\ref{fouzeta}) only in the sign of the last term.
\end{proof}

\medskip

\section{Product of two zeta and related functions } \label{examples}

In this section we evaluate integrals with integrands consisting of products of
two Hurwitz zeta functions. Classical relations for the
Bernoulli polynomials are obtained as corollaries. \\

\no
\begin{Thm}
\label{thm-zetazeta}
Let $z, \, z' \in \npR$.
Then
\begin{align}
\ione \zeta(z',q) \zeta(z,q) dq & = \frac{2 \Gamma(1-z) \Gamma(1-z')}
{(2 \pi)^{2 - z - z'}} \zeta(2 - z - z')
\cos \left( \frac{\pi(z-z')}{2} \right)
\stboth
\label{zzpri1} \\
& = - \zeta(z+z'-1) B(1-z, 1-z') \;
\frac{ \cos(\pi(z-z')/2)}{\cos(\pi(z+z')/2)}.
\stboth
\label{zprima}
\end{align}
\no
Similarly,
\begin{align}
\ione \zeta(z',q) \zeta(z,1-q) dq & = -\frac{2 \Gamma(1-z) \Gamma(1-z')}
{(2 \pi)^{2 - z - z'}} \zeta(2 - z - z')
\cos \left( \frac{\pi(z+z')}{2} \right)
\stboth
\label{zzpri1refl} \\
& = \zeta(z+z'-1) B(1-z, 1-z').
\stboth
\label{zprimarefl}
\end{align}
\end{Thm}
\no
\begin{proof}
The expansion (\ref{fouzeta}) shows that the coefficients
of $\zeta(z',q)$  are given by
\ba
a_{n} & = & \frac{2 \Gamma(1-z')\sin(\pi z'/2)}{(2 \pi)^{1-z'}}
\frac{1}{n^{1-z'}},
\stboth
\nn \\
b_{n} & = & \frac{2 \Gamma(1-z')\cos(\pi z'/2)}{(2 \pi)^{1-z'}}
\frac{1}{n^{1-z'}}.
\stboth
\nn
\ea
\no
Theorem \ref{main} then yields (\ref{zzpri1}). Now use
Riemann's relation (\ref{Riemann1}) for the $\zeta$-function  to obtain
(\ref{zprima}). The proofs of (\ref{zzpri1refl}) and (\ref{zprimarefl}) are
similar. \\
\end{proof}

\medskip

\no
\begin{Example}\label{ex-zetasq}
Let $z \in \npR$. Then  \\
\begin{align}
\ione \zeta^{2}(z,q) dq & =
2 \Gamma^{2}(1-z) ( 2 \pi)^{2z-2} \zeta(2 - 2z)
\stboth
\label{zetasq}\\
\intertext{and}
\ione \zeta(z,q) \zeta(z,1-q) dq & = - 2 \Gamma^{2}(1-z) (2 \pi)^{2z-2}
\zeta(2-2z) \cos(\pi z).
\label{zetasqrefl}
\stboth
\end{align}
\end{Example}
\begin{proof}
Let $z = z'$ in (\ref{zzpri1}) and (\ref{zzpri1refl}).\\
\end{proof}

\no
\begin{Example}
Let $m \in \nnN$. Then
\ba
\ione B_{m}^{2}(q) dq & = & \frac{\abs{B_{2m}} }{ \binom{2m}{m}}.
\label{bsq}
\stboth
\ea
\end{Example}
\begin{proof}
For $m \geq 1$ let $z = 1 - m$ in (\ref{zetasq}).  The case $m=0$ is direct.
\end{proof}
\medskip

\no
\begin{Example}
Let  $m \in \allN$. Then
\ba
\ione \zeta^{2}( -m + \tfrac{1}{2}, q) dq & = &
\left( \frac{(2m)!}{2^{2m} \, m! }
\right)^{2} \; \frac{\zeta(2m+1)}{( 2 \pi)^{2m} }.
\stboth
\label{zetahalf}
\ea
\end{Example}
\begin{proof}
Let $z =  -m + \tfrac{1}{2}$ in (\ref{zetasq}) and use
\ba
\Gamma \left(m + \tfrac{1}{2} \right) & = &
\frac{\sqrt{\pi} (2m)!}{2^{2m} \, m!}.
\stboth
\nn
\ea
\end{proof}

\no
In particular, for $z = -\tfrac{1}{2} \, (m=1)$ we obtain
\ba
\ione \zeta^{2}(-\tfrac{1}{2},q) dq & = & \frac{\zeta(3)}{16 \pi^{2}}.
\stboth
\label{zetasq2}
\ea

\no
{\bf Note}. The integral
\ba
\ione \zeta^{2}( -m + \tfrac{1}{2}, q) dq  \nn
\stboth
\ea
\no
is a rational multiple of $\zeta(2m+1)/\pi^{2m}$. \\

\medskip

The next two examples present special cases of (\ref{zprima}) that
involve integrals of Bernoulli polynomials.  \\

\no
\begin{Example}\label{ex-berzeta}
Let $z \in \npR$ and $m \in \mathbb{N}$. Then \\
\ba
\ione B_{m}(q) \zeta(z,q) dq & = &
(-1)^{m+1} \frac{ m! \, \zeta(z-m) }{(1-z)_{m}},
\stboth
\label{intber1}
\ea
\no
where $(z)_{k} := z(z+1)(z+2) \cdots (z+k-1)$ is the Pochhammer symbol.
\end{Example}

\begin{proof}
Let $z' = 1-m$ in (\ref{zprima}) to produce
\ba
\ione B_{m}(q) \zeta(z,q) dq & = & (-1)^{m+1} \, m \, B(1-z,m) \zeta(z-m).
\stboth \label{exp1}
\ea
The result then follows from $B(1-z,m) = (m-1)!/(1-z)_{m}$.
\end{proof}

\medskip

The next formula appears as $2.4.2.2$ in \cite{pr}. \\

\no
\begin{Example}
Let $n, m \in \mathbb{N}$. Then
\begin{equation}
\begin{split}
\ione B_{m}(q) B_{n}(q) dq & =
\begin{cases}
(-1)^{m+1} \binom{m+n}{m}^{-1} \; B_{m+n} &
\text{ if } m+n \text{ is even,}  \\
0 & \text{ if } m+n \text{ is odd.}
\end{cases}
\label{orthogonality}
\stboth
\end{split}
\end{equation}
\no
The case $m=n$ confirms (\ref{bsq}).
\end{Example}
\begin{proof}
Let $z = 1-n \in - \nnN$ in
(\ref{intber1}) to obtain
\ba
\ione B_{m}(q) B_{n}(q) dq & = &   \frac{(-1)^{m} n \, m! \zeta(1-n-m)}{(n)_{m}}
\stboth
\nn \\
 & = &  \frac{(-1)^{m} m! n! \zeta(n+m)}{(2 \pi)^{n+m} } \;
\; \; 2 \cos \left(\frac{\pi(m+n)}{2} \right)
\stboth
\nn
\ea
\no
using (\ref{Riemann2}).  The vanishing for $n+m$ odd is clear, and for
$n+m$ even the result follows from (\ref{berzeta}). \\
\end{proof}

\no
{\bf Note}. We can write (\ref{orthogonality}) more simply as
\ba
\ione B_{m}(q) B_{n}(q) dq & = &
(-1)^{m+1} \frac{B_{m+n}}{\binom{m+n}{m}},  \nn
\stboth
\ea
\no
recalling that $B_{k} = 0$ for odd $k>1$. \\

We now establish a formula for the moments of $\zeta(z,q)$. \\

\begin{Thm}
\label{thm-mom}
The moments of the Hurwitz zeta function are given
by
\ba
\ione q^{n} \zeta(z,q) dq & = & - n! \sum_{j=1}^{n}
\frac{\zeta(z-j)}{(z-j)_{j} \, (n-j+1)! }
\stboth
\label{mom} \\
& = & n! \sum_{j=1}^{n}  (-1)^{j+1}
\frac{\zeta(z-j)}{(1-z)_{j} \, (n-j+1)!}.
\stboth
\nn
\ea
\end{Thm}
\begin{proof}
We prove (\ref{mom}) by induction.  The
case $n=1$ follows from (\ref{intber1}) and the vanishing of the integral
of $\zeta(z,q)$. For $n > 1$, integration  by parts yields
\ba
\ione q^{n+1} \zeta(z,q) dq & = & \frac{1}{1-z} \ione q^{n+1}
\frac{\partial}{\partial q} \zeta(z-1,q) dq
\stboth
\nn  \\
 & = & \frac{\zeta(z-1)}{1-z} + \frac{(n+1)!}{1-z}
\sum_{k=2}^{n+1} \frac{\zeta(z-k)}{(z-k)_{k-1} \, (n-k+2)!},
\stboth
\nn
\ea
\no
where we have used (\ref{mom}) for power $n$. The final form is obtained
from the identity  $(1-z) \times (z-k)_{k-1}  = -(z-k)_{k}$. \\
\end{proof}

\no
A direct proof of (\ref{mom}) can be given using the expansion of $q^{n}$
in terms of Bernoulli polynomials given in (\ref{polyinv}) and the
evaluation (\ref{intber1}):

\ba
\ione q^{n} \zeta(z,q) dq & = & \frac{1}{n+1}
\sum_{j=0}^{n} \binom{n+1}{j} \ione B_{j}(q) \zeta(z,q) dq
\stboth
\nn \\
 & = & \frac{1}{n+1} \sum_{j=1}^{n} \binom{n+1}{j} (-1)^{j+1}
\frac{j! \, \zeta(z-j)}{(1-z)_{j}}.
\stboth
\nn
\ea

\medskip

Noting the similitude between (\ref{polyinv}) and
(\ref{polyinv1}), the proof above can be imitated to give
\begin{Example}
For $n\in\allN$,
\ba
\label{mom-refl}
\ione (1-q)^{n} \zeta(z,q) dq & = & -n! \sum_{j=1}^{n}
\frac{\zeta(z-j)}{(1-z)_{j} \, (n-j+1)!}.
\stoli
\ea
\end{Example}

\medskip

\no
The special case $z \in -\nnN$ in Theorem \ref{thm-mom} yields the
moments of the Bernoulli polynomials. \\

\no
\begin{Example}
Let $n, m \in \mathbb{N}$. Then
\begin{align}\label{qbernoulli}
\begin{split}
\ione q^{n} B_{m}(q) dq & =
\frac{1}{n+1} \sum_{j=1}^{n} (-1)^{j+1} \frac{\binom{n+1}{j} }
{\binom{m+j}{j}} \; B_{m+j} \stboth \\
& = \frac{n! \, m!}{(n+1+m)!} \sum_{j=1}^{n} (-1)^{j+1} \binom{n+1+m}{n+1-j}
B_{m+j}.
\stboth
\end{split}
\end{align}
\no
\end{Example}
\begin{proof}
Apply (\ref{bernoulli}) to write
\ba
\ione q^{n} B_{m}(q) dq  & = & -m \ione q^{n} \zeta(1-m,q) dq
\stboth
 \nn \\
 & = &
m \sum_{j=1}^{n} (-1)^{j} \frac{\zeta(1 - m - j) (j-1)!}{(m)_{j}} \binom{n}{j-1}
\stboth \nn \\
 & = & (-1)^{n+1} \sum_{j=1}^{n} (-1)^{j+1}
\frac{B_{m+j}}{\frac{m+j}{m} \frac{(m)_{j}}{j!}} \binom{n+1}{j},
\stboth
\nn
\ea
\no
using (\ref{berzeta}) to go from the second to the
third line. The final form follows from the identity
\ba
\frac{m+j}{m} \times \frac{(m)_{j}}{j!} & = & \binom{m+j}{j}. \nn
\stboth
\ea
\end{proof}

\no
{\bf Note}. The results (\ref{intber1}) and (\ref{mom}) are special cases of
the indefinite integrals
\ba
\int B_{m}(q) \zeta(z,q) dq & = & m! \, \sum_{k=1}^{m+1} (-1)^{k+1}
\frac{B_{m+1-k}(q) \zeta(z-k,q)}{(1-z)_{k} (m+1-k)!} \nn \\
& & \nn \\
\int q^{n} \zeta(z,q) dq & = & n! \, \sum_{k=1}^{n+1} (-1)^{k+1}
\frac{q^{n+1-k} \zeta(z-k,q) }{(1-z)_{k} (n+1-k)!} \nn
\ea
\no
discussed in \cite{bem}. \\

\section{The exponential function}

In this section we evaluate the Hurwitz transform of the exponential function.
The result is expressed in terms of the transcendental function
\ba
F(x,z) & := & \sum_{n=0}^{\infty} \zeta(n+2-z) x^{n}, \quad \text{ for }
|x| < 1.
\stboth
\label{funcionF}
\ea

\medskip

\no
\begin{Example}
Let $ z \in \npR$ and $\abs{t} < 1$. Then
\ba
\quad \quad \ione e^{2 \pi tq} \zeta(z,q) dq & = & 2 (1 - e^{2 \pi t}) \;
\frac{\Gamma(1-z) }{(2 \pi)^{2-z}} \, \times
\realpart\left[ e^{\pi i z/2} F(it,z) \right],
\stboth
\label{exponen}
\ea
\no
where $F(x,z)$ is given in (\ref{funcionF}). \\
\end{Example}

\begin{proof}
The generating function for the Bernoulli polynomials
(\ref{gener}) yields
\ba
e^{qt} & =& \frac{e^{t} -1}{t} \sum_{n=0}^{\infty} B_{n}(q)
\frac{t^{n}}{n!},
\stboth
\nn
\ea
\no
so that
\ba
\ione e^{qt} \zeta(z,q) dq & = &
\frac{e^{t} -1}{t} \sum_{n=0}^{\infty} \frac{t^{n}}{n!}
\ione B_{n}(q) \zeta(z,q) dq. \nn
\stboth
\ea
\no
Since $B_{0}(q)=1$ and $\zeta(z,q)$ integrates to $0$, the above sum
effectively starts at $n=1$. Thus (\ref{exp1}) gives
\ba
\ione e^{qt} \zeta(z,q) dq & = &  (e^{t}-1) \sum_{n=0}^{\infty}
(-1)^{n} \frac{t^{n}}{n!} B(1-z,n+1) \zeta(z-n-1), \nn
\stboth
\ea
\no
which can be written as
\begin{equation}
\ione e^{qt} \zeta(z,q) dq  =  \frac{2(e^{t}-1) \Gamma(1-z)}
{(2 \pi)^{2-z}}  \sum_{n=0}^{\infty}  (-1)^{n} \!\left( \frac{t}{2 \pi}
\right)^{n}  \zeta(n+2-z) \cos \left(\frac{\pi(z-n)}{2} \right) \nn
\stboth
\end{equation}
\no
using (\ref{Riemann2}) and (\ref{betagamma}). Now
replace $t$ by $2 \pi t$ and use the evaluation $\cos(\pi(z-n)/2) =
\realpart(e^{i \pi (z-n)/2})=\realpart((-i)^n e^{i \pi z/2})$ to
yield the final result. \\
\end{proof}

The next example results from $z \in - \nnN$ in (\ref{exponen}). It
appears in \cite{pr}:
$2.4.1.4$.\footnote{The factor $m!$ in (\ref{expber}) is missing in \cite{pr}.} \\

\no
\begin{Example}
Let $m \in \allN$ and $\abs{t}<1$. Then
\begin{multline}\label{expber}
\ione e^{2 \pi t q} B_{m}(q) dq  =  \frac{(-1)^{m}(e^{2 \pi t}-1) \, m!}
{(2 \pi t)^{m+1}} \times \\
\times \left[ 1  - \pi t \;
\text{ coth} (\pi t) - 2 \sum_{r=1}^{\lfloor{ \tfrac{m}{2} \rfloor} }
(-1)^{r} \zeta(2r) t^{2r} \right].
\stboth
\end{multline}
\no
\end{Example}
\begin{proof}
We discuss the case $m = 2k+1$; the case of $m$ even is similar. Let
$z = 1-m = -2k$
in  (\ref{exponen}). Then
\ba
\realpart\left[ e^{\pi i z/2} F(it, z) \right]  & = &
(-1)^{k} \realpart \left[ F(it, -2k) \right]
\stboth
 \nn \\
& = & (-1)^{k} \sum_{r=0}^{\infty} \zeta( 2r+2+2k) (-1)^{r} t^{2r}
\stboth
\nn \\
 & & \nn \\
 & = & - t^{2k+2} \left[ \frac{1}{2} - \frac{\pi t}{2} \text{ coth } \pi t -
\sum_{r=1}^{k}
(-1)^{r} \zeta(2r) t^{2r} \right], \nn
\stboth
\ea
\no
where we have employed the identity
\ba
\text{ coth } \pi x & = & \frac{1}{\pi x} - \frac{2}{\pi x}
\sum_{r=1}^{\infty} (-1)^{r} \zeta(2r) x^{2r} \label{coth}
\stboth
\ea
\no
that appears in \cite{atlas}, 3:14:5.  \\
\end{proof}

\no
\section{The logsine function}

This section contains examples involving the function $\ln (\sin \pi q)$.
The standard tables \cite{gr} and \cite{pr} contain very few examples of this
type. See sections $4.224$ and $4.322$. Some of the evaluations presented
here are computable by
Mathematica 4.0. \\

\no
\begin{Example}
Let $z \in \npR$. Then
\ba
\ione \ln ( \sin \pi q ) \, \zeta(z,q) dq & = &
- \frac{\Gamma(1-z)}{(2 \pi)^{1-z}} \sin \left( \frac{\pi z}{2} \right)
\zeta(2-z)
\label{lotszeta1}
\stboth \\
& =  & - \frac{\zeta(z) \, \zeta(2-z)}{2 \zeta(1-z)},
\stboth
\label{lotszeta}
\ea
\no
where the second result follows from (\ref{lotszeta1}) when $ z \neq 0$
by use of (\ref{Riemann1}).
\end{Example}
\begin{proof}
The Fourier coefficients of $\ln ( \sin \pi q)$ are
\ba
\ione \ln ( \sin \pi q ) \sin( 2 n \pi q) dq  & = & 0
\stboth
\nn
\ea
\no
and
\ba
\ione \ln ( \sin \pi q ) \cos( 2 n \pi q) dq  & = &
\begin{cases}
- \ln 2& \text{ if } n = 0, \\
- \tfrac{1}{2n}& \text{ if } n > 0.
\end{cases}
\stboth
\nn
\ea
\no
These appear in \cite{gr} 4.384. Thus (\ref{lotszeta1}) follows from
Theorem \ref{main}.  \\
\end{proof}

\medskip

\no
\begin{Example}
Let $m \in \mathbb{N}$. Then
\begin{equation}\label{logsine}
\begin{split}
\ione \ln ( \sin \pi q ) \, B_{m}(q) dq & =
\begin{cases}
(-1)^{m/2}  (2 \pi)^{-m} m! \zeta(m+1) & \quad \text{ if } m
\text{ is even,} \\
0 & \quad \text{ if } m \text{ is odd}.
\end{cases}
\stboth
\end{split}
\end{equation}
\end{Example}
\begin{proof}
Let $z = 1-m \in - \nnN$ in
(\ref{lotszeta}) giving
\ba
\ione \ln ( \sin \pi q ) \, B_{m}(q) dq & = &
\frac{m \zeta(1-m) \zeta(1+m)}{2 \zeta(m)}.
\stboth
\label{zetam}
\ea
\no
Now use (\ref{Riemann2}) to obtain the result. \\
\end{proof}

\no
{\bf Note}. The integral
\ba
\ione \ln ( \sin \pi q ) \, B_{2m}(q) dq
\stboth
\ea
\no
is a rational multiple of $\zeta(2m+1)/\pi^{2m}$.  \\

\medskip

The next example evaluates the moments of $\ln( \sin \pi q)$.  \\

\no
\begin{Example}
Let $n \in \nnN$. Then
\ba
\ione q^{n} \ln (\sin \pi q) dq & = &
- \frac{\ln 2}{n+1} + n! \sum_{k=1}^{\lfloor{ \tfrac{n}{2} \rfloor} }
\frac{(-1)^{k} \zeta(2k+1)}{(2 \pi)^{2k} \, (n+1 - 2k)!}.
\stboth
\label{momlogsine}
\ea
\end{Example}
\begin{proof}
Using (\ref{polyinv}) we have
\ba
\ione q^{n} \ln ( \sin \pi q) dq & = &
\frac{1}{n+1} \sum_{j=0}^{n} \binom{n+1}{j}
\ione \ln ( \sin \pi q ) \, B_{j}(q) dq.
\stboth
 \nn
\ea
\no
The result now follows by (\ref{logsine}) and
the classical value
\ba
\ione \ln (\sin \pi q )dq & = & - \ln 2.
\stboth
\label{classic}
\ea
An elementary evaluation of (\ref{classic}) appears in
\cite{arora}.  \\
\end{proof}

\medskip

\no
{\bf Note}. The integral
\ba
\ione q^{n} \ln ( \sin \pi q ) dq
\stboth
\ea
\no
is a rational linear combination of $\ln 2$ and
$\{ \zeta(2k+1)/\pi^{2k}: 1 \leq k \leq \lfloor{\frac{n}{2} \rfloor} \}$.  \\

\medskip

\no
The first few cases are

\ba
\ione q \, \ln ( \sin \pi q ) dq & = &   - \frac{1}{2} \ln 2,
\stboth
\label{example4}  \\
\ione q^{2}  \, \ln ( \sin \pi q ) dq & = &  -\frac{1}{3} \ln 2 -
\frac{\zeta(3)}{2 \pi^{2}}, \stboth \nn \\
\ione q^{3}  \, \ln ( \sin \pi q ) dq & = &  - \frac{1}{4} \ln 2 -
\frac{3 \zeta(3)}{4 \pi^{2}}, \stboth \nn \\
\ione q^{4} \ln ( \sin \pi q) dq & = &
- \frac{1}{5} \ln 2 - \frac{\zeta(3)}{\pi^{2}} + \frac{3 \zeta(5)}{2 \pi^{4}}.
\stboth \nn
\ea
\no
These evaluations can be confirmed by Mathematica 4.0. \\

\no
\section{The loggamma function}
\label{sec-loggamma}

This section contains evaluations involving the function
$\ln \Gamma(q)$. None of the examples presented here
were computable by a symbolic language.  \\

\no
\begin{Example}
Let $z \in \npR$. Then
\ba
\;\;\;\; \ione \ln \Gamma(q) \, \zeta(z,q) dq & = &
\frac{\Gamma(1-z)}{(2 \pi)^{2-z}} \zeta(2-z)
\stboth
\label{loggamma}  \\
    & &\times  \left[ \pi \, \sin \left( \frac{\pi z}{2} \right)  +
2 \cos \left( \frac{\pi z}{2} \right)
\left\{ A  -  \frac{\zeta'(2-z)}{\zeta(2-z)} \right\} \right], \nn
\ea
\no
where
\ba
A & := & 2 \lp + \gamma =
-2 \frac{d}{dz} \left( \zeta(z) \Gamma(1-z) \right) \Big{|}_{z=0}
\stboth
\label{letterA}
\ea
\no
and $\gamma$ is Euler's constant.
\end{Example}

\begin{proof}
The Fourier coefficients of $\ln \Gamma(q)$ appear in \cite{gr} $6.443.1$ and
$6.443.3$ as
\ba
\ione \ln \Gamma(q) \, \sin(2 \pi n q) dq & = & \frac{ A +
\ln n }{2 \pi n},\quad n\in\allN,
\stboth
\label{logga1} \\
\ione \ln \Gamma(q) \, \cos(2 \pi n q) dq & = & \frac{1}{4n},\quad n\in\allN.
\stboth
\label{logga2}
\ea
\no
Thus
\ba
a_{n} = \frac{1}{2n} & \text{ and } & b_{n} =
\frac{A + \ln n }{ \pi n},
\stboth
\nn
\ea
\no
where $A$ is defined in (\ref{letterA}).  The evaluations
\ba
\sum_{n=1}^{\infty} \frac{1}{2n^{2-z}}  =
\frac{1}{2} \zeta(2-z)  & \text{ and } &
\sum_{n=1}^{\infty} \frac{ A + \ln n}{n^{2-z}} =
A \zeta(2-z) - \zeta'(2-z)
\stboth
\nn
\ea
\no
yield (\ref{loggamma}).
\end{proof}

\medskip

\no
{\bf Note}. The integral
\ba
\ione \ln \Gamma(q) \, \cos((2n+1) \pi q) dq & = &
\frac{2}{\pi^{2}} \left( \frac{\gamma + 2 \lp }{(2n+1)^{2}} +
2 \sum_{k=2}^{\infty} \frac{\ln k}{4k^{2} -(2n+1)^{2}} \right), \nn
\ea
\no
a companion to (\ref{logga2}), was evaluated
by K\"{o}lbig in \cite{kolbig3}. This
was recorded as $0$ as late as in the fourth edition of
\cite{gr}. The fifth edition contains the correct value. \\
\medskip

\no
\begin{Example}\label{Ex-bergamma}
Let $ m \in \mathbb{N}$. Then
\begin{align}
\ione B_{2m}(q) \ln \Gamma(q) dq & = (-1)^{m+1} \frac{(2m)! \zeta(2m+1)}
{2(2 \pi)^{2m}} = - \zeta'(-2m), \stboth \label{special1} \\
\ione B_{2m-1}(q) \ln \Gamma(q) dq & =
\frac{B_{2m}}{2m} \times
\left[ \frac{\zeta'(2m)}{\zeta(2m)} - A  \right].  \stboth \label{special2}
\end{align}
\end{Example}
\begin{proof}
\no
Replace in (\ref{loggamma}) the variable $z$ by $1-2m$ and $2-2m$
respectively. Then use (\ref{berzeta2}) in the first case and
(\ref{berzeta}) in the second.
\\
\end{proof}

\no
{\bf An alternative approach}. The
evaluation in Example \ref{Ex-bergamma} can also be obtained by integrating
\ba
\left. {\frac{d}{{dz}} \zeta(z,q)} \right|_{z=0} & = & \ln \Gamma (q) -
\lp
\stboth
\label{loggamma0}
\ea
\no
to produce
\begin{align}
\ione B_{m}(q) \ln \Gamma(q) dq & = & \lp \ione B_{m}(q)dq +
\frac{d}{dz} \Bigg{|}_{z=0} \ione B_{m}(q) \zeta(z,q) dq. \label{log0}
\stboth
\end{align}
\no
The relation (\ref{loggamma0}) can be found in \cite{gr} $9.533.3$.  To
evaluate (\ref{log0}) differentiate  (\ref{intber1}) to produce
\ba
\ione B_{m}(q) \ln \Gamma(q) dq & = & \lp \delta_{m,0} +
(1 - \delta_{m,0}) (-1)^{m+1} \left[ H_{m} \zeta(-m) + \zeta'(-m) \right].
\nn
\stboth
\ea
\no
Here $\delta_{m,0}$ is Kronecker's delta and $H_{m} = 1 + \tfrac{1}{2} +
\cdots + \tfrac{1}{m}$ is the $m$-th harmonic number.
Use has been made of the result
\ba
\left. {\frac{d}{{dz}}(1 - z)_k } \right|_{z = 0}  =  - k!H_k.
\label{pochderiv}
\stboth
\ea
According to the parity of $m$ we have
\begin{equation}
\begin{split}
\ione B_{m}(q) \ln \Gamma(q) dq & =
\begin{cases}
- \zeta'(-m) & m = 0, \,2, \,4, \cdots, \\
H_{m} \zeta(-m) + \zeta'(-m)  & m = 1, \, 3, \cdots
\end{cases}
\label{special3}
\stboth
\end{split}
\end{equation}
(for $m=0$ we have used (\ref{berzeta3})). The
result (\ref{special3}) for odd $m$ is seen to be equivalent
to (\ref{special2}) after use of the identity
\ba\label{Riemann3}
\frac{{\zeta '(1 - 2k)}}{{\zeta (1 - 2k)}} + \frac{{\zeta
'(2k)}}{{\zeta (2k)}} = \ln 2\pi  + \gamma  - H_{2k - 1}, \quad
k\in\allN,
\stboth
\ea
which can be derived by differentiating Riemann's relation
(\ref{Riemann1}) and evaluating at $s=2k$.
\medskip

\no
\begin{Example}
The case $m=1$ in (\ref{special2}) yields, using
$\zeta(2) = \pi^{2}/6$,
\ba
\ione ( q - \tfrac{1}{2} ) \, \ln \Gamma(q) dq & = & \frac{1}{12}
\left( \frac{6 \zeta'(2) }{\pi^{2}} - 2 \lp - \gamma \right).
\stboth
\label{zzero}
\ea
\no
The case $m=1$ in (\ref{special1}) gives
\ba
\ione ( q^2 - q + \tfrac{1}{6}) \; \ln \Gamma(q)  dq & = &
 \frac{\zeta(3)}{4 \pi^{2}}.
\stboth \label{zeta3}
\ea
\end{Example}

\medskip

\no
\begin{Example}
Let $n \in \mathbb{N}$. Then
\begin{align}\label{loggamma1}
\begin{split}
\ione q^{n} \ln \Gamma(q) dq & =
\frac{1}{n+1} \sum_{k=1}^{\lfloor{ \tfrac{n+1}{2} \rfloor}  }
(-1)^{k} \binom{n+1}{2k-1}  \frac{(2k)!}{k (2 \pi)^{2k}}
\left[ A  \zeta(2k) - \zeta'(2k) \right] \\
&-\frac{1}{n+1} \sum_{k=1}^{\lfloor{ \tfrac{n}{2} \rfloor} }
(-1)^{k} \binom{n+1}{2k}  \frac{(2k)!}{2 (2 \pi)^{2k}}
\zeta(2k+1) +  \frac{\lp}{n+1}.
\stboth
\end{split}
\end{align}
\end{Example}

\begin{proof}
Use the expression (\ref{polyinv}) to write
\ba
\ione q^{n} \ln \Gamma(q) dq  & = & \frac{1}{n+1}
\sum_{j=1}^{n} \binom{n+1}{j} \ione \ln \Gamma(q) B_{j}(q) dq
\stboth
 \nn  \\
   & + & \frac{1}{n+1} \ione \ln \Gamma(q) dq. \nn
\ea
\no
The value
\ba
\ione \ln \Gamma(q) dq & = &  \lp \label{nice}
\stboth
\ea
\no
is then obtained from (\ref{loggamma0}) and (\ref{vanishing}).
The result  now follows from (\ref{special1}) and
(\ref{special2}).
\end{proof}

The formula (\ref{loggamma1}) yields \\
\ba
\ione q \ln \Gamma(q) dq & = & \frac{\zeta'(2)}{2 \pi^{2}} +
\frac{1}{3} \lp - \frac{\gamma}{12},
\stboth
\nn \\
& & \nn \\
\ione q^{2} \ln \Gamma(q) dq & = & \frac{\zeta'(2)}{2 \pi^{2}} +
\frac{\zeta(3)}{4 \pi^{2}} +
\frac{1}{6} \lp - \frac{\gamma}{12},
\stboth
\nn \\
& & \nn \\
\ione q^{3} \ln \Gamma(q) dq & = & \frac{\zeta'(2)}{2 \pi^{2}} +
\frac{3 \zeta(3)}{8 \pi^{2}}  - \frac{3 \zeta'(4)}{4 \pi^{4}} +
\frac{1}{10} \lp - \frac{3 \gamma}{40}.
\stboth
\nn
\ea
\no
None of these examples could be evaluated symbolically.  \\

\medskip

\no
{\bf Note}. The  integral
\ba
\ione q^{n} \ln \Gamma(q) dq  \nn
\stboth
\ea
\no
is a rational linear combination of
\ba
\left\{ \gamma, \lp, \;
\frac{\zeta'(2k)}{\pi^{2k}} \; 1 \leq k \leq \lfloor{\frac{n}{2}\rfloor},
\;  \frac{\zeta(2k+1)}{\pi^{2k}},
\; 1 \leq k \leq  \lfloor{\frac{n+1}{2} \rfloor} \right\}. \nn
\stboth
\ea

\medskip

\no
{\bf Note}. Gosper \cite{gosper} presents a series of interesting
evaluations of definite integrals of $\ln \Gamma(q)$. For example
\ba\label{gosper1}
\int_{0}^{1/2} \ln \Gamma(q+1) dq & = & \frac{\gamma}{8}  +
\frac{3 \lp }{4}
- \frac{13 \ln 2}{24} - \frac{3 \zeta'(2)}{4 \pi^{2}} - \frac{1}{2}
\st
\ea
\no
and
\ba
\label{gosper2}
\int_{0}^{1/4} \ln \Gamma(q+1) dq & = & \frac{3\gamma}{32}  +
\frac{7 \lp}{16} - \frac{\ln 2}{2} - \frac{9 \zeta'(2)}{16 \pi^{2}} +
\frac{G}{4 \pi} - \frac{1}{4},
\st
\ea
\no
where
\ba
G & := & \sum_{n=0}^{\infty} \frac{(-1)^{n}}{(2n+1)^{2}}  \label{catalan}
\stboth
\ea
\no
is Catalan's constant. See Section \ref{inthalf} for an alternative proof
of (\ref{gosper1}). \\

\no
{\bf Note}. The results discussed here are special cases of the indefinite
integral
\ba
\int q^{n} \ln \Gamma(q) dq & = & -\zeta'(0) \frac{q^{n+1}}{n+1} \nn \\
          & + & n! \sum_{k=1}^{n+1} (-1)^{k+1} \frac{q^{n+1-k}}{k! (n+1-k)!}
\left[ \zeta_{z}(-k,q) - \frac{H_{k}}{k+1} B_{k+1}(q) \right], \nn
\stoli
\ea
\no
where $H_{k}$ is the $k$-th harmonic number and
\ba
\zeta_{z}(-k,q) & := & \left. {\frac{\partial}{{\partial z}}}\right|_{z=-k}
\zeta(z,q).
\stboth
\ea
\no
These results can be expressed in terms Gosper's {\em negapolygammas}
$\psi_{-k}(q)$  \cite{gosper} in view of the relation
\ba
\zeta_{z}(-k,q) & = & \frac{H_{k}}{k+1} B_{k+1}(q) + q^{k} \zeta'(0)
+k!\,
\psi_{-k}(q), \nn
\stoli
\ea
\no
where
\ba
\psi_{-1}(q) & = & \ln \Gamma(q), \nn \\
\psi_{-k}(q) & = & \int \psi_{-k+1}(q) dq,\quad k\ge 2. \nn
\stboth
\ea
Details will appear in \cite{bem}. \\

\no
\section{Differentiation results}

In this section we discuss evaluation of certain integrals that appear
from (\ref{zzpri1}) after differentiation with respect to the parameters
$z$ and $z'$.  The special values $z=0$ and $z'=0$ produce evaluations
containing the loggamma function, in view of (\ref{loggamma0}).
In particular, as was pointed out earlier, the
result
\ba
\ione \ln \Gamma(q) dq & = & \lp \label{logpi}
\stboth
\ea
\no
follows directly from (\ref{loggamma0}).  The integrals considered
here complement those considered in Section \ref{sec-loggamma}. \\

\begin{Prop}
For $z, z' \in \npR$ we have
\begin{multline}\label{derivmain}
\ione \frac{d}{dz} \zeta(z,q) \zeta(z',q) dq = -
\frac{2 \Gamma(1-z) \Gamma(1-z')}{(2 \pi)^{2-z-z'}} \zeta(2-z-z') \, \cos \omega
\\
\times \left[ \frac{\zeta'(2-z-z')}{\zeta(2-z-z')}
 + \frac{\pi}{2} \tan \omega
- 2 \lp + \psi(1-z) \right],
\stboth
\end{multline}
\no
where $\omega = \pi(z-z')/2$ and $\psi(z)$ is the digamma
function defined in (\ref{digamma}).
\end{Prop}
\begin{proof}
Direct differentiation of (\ref{zzpri1}). \\
\end{proof}

\no
In particular, for $z=z'=0$ we obtain  (\ref{zzero}).  \\

\no
\begin{Example}
Differentiating (\ref{zzpri1}) with respect to $z$ and
then $z'$, evaluating at $z = z'=0$, and using (\ref{loggamma0})
yields
\\
\begin{multline}\label{loggammasquared}
\int_0^1 \left( {\ln \Gamma (q)} \right)^2 dq = \frac{\gamma ^2
}{12}  + \frac{\pi ^2}{48} + \frac{1}{3}\gamma \lp
+ \frac{4}{3}\left( {\lp } \right)^2 \\
- (\gamma  + 2 \lp )\frac{\zeta '(2)}{\pi ^2} +
\frac{\zeta ''(2)}{2\pi ^2 }.
\stboth
\end{multline}
\end{Example}

\bigskip

\no
\begin{Example}
Differentiating (\ref{lotszeta1}) with respect
to $z$ and then setting $z=0$ yields, after using
(\ref{classic}),
\ba
\int_0^1 {\ln \sin \pi q\ln \Gamma (q)\,dq = }  - \ln
2 \lp  - \frac{{\pi ^2 }}{{24}}.
\stboth
\label{logslogg}
\ea
\end{Example}

\medskip

\no
\section{An expression for Catalan's constant}

In his discussion of Entry 17(v) of Chapter 8 of
Ramanujan's Notebooks,
Berndt \cite{ram1} page 200, introduces the function
\ba
G(z,q)  & := & \zeta(z,q) - \zeta(z,1-q) \label{G}
\ea
\no
and gives its Fourier expansion
\ba
G(z,q) & = & 4 \Gamma(1-z) \cos \left(\frac{\pi z}{2} \right)
\, \sum_{k=1}^{\infty}
\frac{\sin(2 \pi k q)}{(2 \pi k)^{1-z}}. \label{bruceexp}
\ea
\no
This is an immediate consequence of the Fourier
expansion (\ref{fouzeta}) for $\zeta(z,q)$. \\

In terms of $G(z,q)$ we can define an anti-symmetrized Hurwitz
transform,
\ba
\mathfrak{H}_A(f) & :=  & \frac{1}{2} \ione  f(w, q) G(z,q) dq.
\stboth
\label{antihurtrans1}
\ea
\no
It is straightforward to show that for a function $f(w,q)$ with
Fourier expansion as in Theorem (\ref{main}) one obtains
\ba
\frac{1}{2} \ione  f(w,q) G(z,q) dq =
\frac{\Gamma(1-z)\cos \left(\pi z/2 \right) }{(2 \pi)^{1-z}}
\sum_{n=1}^{\infty} \frac{b_{n}(w)}{n^{1-z}}.
\label{antihurtrans2}
\ea
As a particular example we compute the anti-symmetrized Hurwitz
transform of $\text{sec}(\pi q)$ and obtain as a corollary an
expression for Catalan's constant.  \\

\no
\begin{Example}\label{Ex-Gsec}
The anti-symmetrized Hurwitz transform of $\text{sec}(\pi q)$ is
\ba
\frac{1}{2} \ione \frac{\zeta(z,q) - \zeta(z,1-q)}{\cos(\pi q)} dq & = &
\frac{16 \Gamma(1-z) \cos(\pi z/2)}{(2 \pi)^{2-z}}
\sum_{n=1}^{\infty} \frac{(-1)^{n+1} }{n^{1-z}} \sum_{k=0}^{n-1}
\frac{(-1)^{k}}{2k+1}. \nn
\stboth
\ea
\end{Example}
\begin{proof}
In \cite{gr} $3.612.5$  we find
\ba
\int_{0}^{1} \frac{\sin(2n \pi q)}{\cos(\pi q)} dq & = &
(-1)^{n+1} \frac{4}{\pi} \sum_{k=0}^{n-1} \frac{ (-1)^{k} }{2k+1}.
\label{gr1}
\stboth
\ea
\no
A straightforward application of (\ref{antihurtrans2})
completes the proof. \\
\end{proof}

\medskip

The special case $z=0$ yields the following result.

\begin{Prop}
The Catalan constant $G$, defined in (\ref{catalan}), is given by
\ba
G  & = & \sum_{n=1}^{\infty} \frac{(-1)^{n+1}}{n} \sum_{k=0}^{n-1}
\frac{(-1)^{k}}{2k+1}.
\label{newcata}
\stboth
\ea
\end{Prop}
\begin{proof}
Put $z=0$ in example \ref{Ex-Gsec} to obtain
\ba
\ione \frac{\tfrac{1}{2} - q}{\cos(\pi q)} dq & = &
\frac{4}{\pi^{2}}
\sum_{n=1}^{\infty} \frac{(-1)^{n+1} }{n} \sum_{k=0}^{n-1}
\frac{(-1)^{k}}{2k+1}. \label{intcata}
\stboth
\ea
\no
The change of variable $t  = \pi\left(\tfrac{1}{2} - q\right)$ then produces
\ba
\ione \frac{\tfrac{1}{2} - q}{\cos(\pi q)} dq  =
\frac{2}{\pi^{2}} \int_{0}^{\pi/2} \frac{t}{\sin t} dt
=\frac{4G}{\pi^2} \label{cata1}
\stboth
\ea
\no
since the second integral equals $2G$.
\end{proof}

\no
{\bf Note}. The direct symbolic evaluation of the integral in
(\ref{intcata}) yields
\ba
\ione \frac{\tfrac{1}{2} - q}{\cos(\pi q)} dq & = &
\frac{1}{16 \pi} \left[ \frac{32 G}{\pi} - \; {_{4}F_{3} }\left(
\begin{matrix}
1 \;  1 \;   \tfrac{3}{2} \;  \tfrac{3}{2}  \\
2 \;  2 \;  2
\end{matrix}
; 1 \right) + 16 \ln 2  \right], \nn
\stboth
\ea
\no
and thus
\ba
G & = & \frac{\pi}{32} \left[ 16 \ln 2 -
{_{4}F_{3}} \left(
\begin{matrix}
1 \;  1 \;   \tfrac{3}{2} \;  \tfrac{3}{2}  \\
2 \;  2 \;  2
\end{matrix}
;1 \right) \right]. \label{catalan2}
\stboth
\ea
\no
This form of Catalan's constant appears in \cite{pr}: $7.5.3.120$, and
(\ref{cata1}) is Entry $14$ in the list
of expressions for $G$ compiled by Adamchik in
\cite{adamcata}, but (\ref{newcata}) does not appear
there.
\no
\begin{Example}
Let $m \in \nnN$. Then
\begin{align}
\begin{split}
\ione \text{sec}(\pi q) B_{2m+1}(q) dq  & =
(-1)^{m+1} \frac{16(2m+1)!}{(2 \pi)^{2m+2}}
\sum_{n=1}^{\infty} \frac{(-1)^{n+1}}{n^{2m+1}} \sum_{k=0}^{n-1} \frac{(-1)^{k}}
{2k+1}
\label{secant1}
\stboth
\end{split} \\
& = (-1)^{m} \frac{4(2m+1)!}{(2 \pi)^{2m+2}}  \times
\sum\nolimits^{*}  \frac{\psi(n/2+1/4)} {n^{2m+1}},
\label{secant2}
\stboth
\end{align}
\no
where the sum extends over $n \in \mathbb{Z}, \, n \neq 0$.
\end{Example}
\begin{proof}
The value $z = -2m$ in Example \ref{Ex-Gsec} yields
(\ref{secant1}). To prove (\ref{secant2}) it is enough to
establish the identity
\ba
\sum\limits_{n = 1}^\infty  {\frac{{( - 1)^{n + 1} }}{{n^{2m + 1}
}}\sum\limits_{k = 0}^{n - 1} {\frac{{( - 1)^k }}{{2k + 1}}} }  =
- \frac{1}{4}\sum\nolimits^* {\frac{{\psi (n/2 + 1/4)}}{{n^{2m + 1} }}}.
\label{catalan-m}
\stboth
\ea
\no
The internal sum in (\ref{catalan-m}) can be written as
\ba
\sum\limits_{k = 0}^{n - 1} {\frac{{( - 1)^k }}{{2k + 1}}}  =
\frac{\pi }{4} + \frac{1}{4}( - 1)^n \left[ {\psi (n/2 + 1/4) -
\psi (n/2 + 3/4)} \right].
\stboth
\ea
\no
Logarithmic differentiation of the reflection formula
(\ref{reflec}) for the gamma function yields
\ba
\psi(1-x) & = & \psi(x) + \pi \, \text{cotg} \, \pi x, \nn
\stboth
\ea
\no
so that, evaluating at $x=1/4-n/2$,
\ba
\psi(1/4+n/2) - \psi(3/4+n/2) & = & \psi(1/4+n/2)-\psi(1/4-n/2) -
(-1)^{n} \pi.  \nn
\stboth
\ea
\no
Thus
\ba
\sum\limits_{k = 0}^{n - 1} {\frac{{( - 1)^k }}{{2k + 1}}}  =
\frac{1}{4}( - 1)^n \left[ {\psi (1/4 + n/2) - \psi (1/4 - n/2)}
\right]
\stboth
\ea
and (\ref{catalan-m}) is established. \\
\end{proof}

\medskip

\section{Clausen and related functions}

In this section we evaluate the Hurwitz transform of the Clausen functions
$\cl_{n}(q)$. These functions are defined by
\ba
\cl_{2n}(x) & :=  & \sum_{k=1}^{\infty} \frac{\sin kx }{k^{2n}}, \quad
 n \geq 1
\label{claueven}
\stboth
\ea
\no
and
\ba
\cl_{2n+1}(x) & :=  &
\sum_{k=1}^{\infty} \frac{\cos kx }{k^{2n+1}}, \quad  n \geq 0.
\label{clauodd}
\stboth
\ea
Extensive information about these functions can be found in
\cite{dilog}, chapter 4. For example,
\ba
\cl_{1}(x) & = & - \ln | 2 \sin (x/2) |. \nn
\stboth
\ea
\no
More generally, one can define the Clausen functions in terms
of the polylogarithm on the unit circle as
\ba
\cl_{2n} (x):=& \imagpart\li_{2n}
(e^{ix} ),
\label{claueven1}
\stboth
\\
\cl_{2n + 1} (x):=& \realpart\li_{2n + 1} (e^{ix}),
\nn
\stboth
\ea
where, for $\abs{z}\le 1$,
\ba\label{def-polylog}
\li_n (z): = \sum\limits_{k = 1}^\infty  {\frac{{z^k }}{{k^n
}}},\quad n\in\allN.
\stboth
\ea
The Fourier expansion of $\li_n (z)$ on the unit circle,
\ba\label{fou-polylog}
\li_n (e^{2\pi qi} ) = \sum\limits_{k = 1}^\infty  {\frac{{\cos
(2\pi kq)}}{{k^n }}}  + i\sum\limits_{k = 1}^\infty  {\frac{{\sin
(2\pi kq)}}{{k^n }}},\quad 0\le q < 1,
\stboth
\ea
leads us, in view of Theorem \ref{main}, to the  next example.  \\

\begin{Example}
Let $z \in \npR$. Then
\ba\label{int-polylog}
\ione\li_n (e^{2\pi qi} )\zeta (z,q)dq = \frac{{\Gamma (1 - z)}}{{(2\pi
)^{1 - z} }}e^{i{\textstyle{\pi  \over 2}}(1 - z)} \zeta (1 - z +
n).
\stboth
\ea
\end{Example}
\no
As immediate consequences we have the next three examples. \\

\no
\begin{Example}
Let $z \in \npR$. Then
\ba
\ione \cl_{2n}(2 \pi q) \zeta(z,q) dq & = &
\frac{\Gamma(1-z) \cos(\pi z/2)}{(2 \pi)^{1-z}} \zeta(1 - z + 2n)
\label{clausen1}
\stboth
\ea
\no
and
\ba
\ione \cl_{2n+1}(2 \pi q) \zeta(z,q) dq & = &
\frac{\Gamma(1-z) \sin(\pi z/2)}{(2 \pi)^{1-z}} \zeta(2 - z + 2n).
\label{clausen2}
\stboth
\ea
\end{Example}

\medskip

\no
\begin{Example}
Let $m \in \mathbb{N}$. Then
\begin{align}
\ione B_{m}(q) \cl_{2n}(2 \pi q) dq & =
\begin{cases}
0  & \text{ if } m \text{ is even,} \\
(-1)^{\tfrac{m+1}{2}} \, m! ( 2 \pi)^{-m} \, \zeta(m + 2n) & \text{ if  } m
\text{ is odd,}
\end{cases}
\nn
\stboth
\end{align}
\no
and
\begin{align}
\ione B_{m}(q) \cl_{2n+1}(2 \pi q) dq & =
\begin{cases}
0  & \text{ if } m \text{ is odd,} \\
(-1)^{\tfrac{m}{2}+1} \, m! ( 2 \pi)^{-m} \, \zeta(m + 2n+1) & \text{ if  } m
\text{ is even}.
\end{cases}
\nn
\stboth
\end{align}
\end{Example}

\medskip

\no
\begin{Example}
Let $m \in \mathbb{N}$. Then
\ba
\ione q^{m} \cl_{2n}(2 \pi q) dq & = &
m! \, \sum_{j=0}^{\lfloor{ \tfrac{m-1}{2} \rfloor}} \frac{(-1)^{j+1}
\, \zeta(2n+2j+1) }{(m-2j)! \, (2 \pi)^{2j+1} } \nn
\stboth
\ea
\no
and
\ba
\ione q^{m} \cl_{2n+1}(2 \pi q) dq & = &
m! \, \sum_{j=1}^{\lfloor{ \tfrac{m}{2} \rfloor}} \frac{(-1)^{j+1}
\, \zeta(2n+2j+1) }{(m-2j+1)! \, (2 \pi)^{2j} }. \nn
\stboth
\ea
\end{Example}

\medskip

\no
\section{A function from Berndt's work on Ramanujan Notebooks}

In his reinterpretation of Entry 3, Chapter 9 of
Ramanujan's Notebooks, Berndt \cite{ram1}, page 235, introduces the functions
\ba
S_{N}(x) & := & \sum_{n=0}^{\infty}  \frac{(-1)^{n} \sin(2n+1)x}{(2n+1)^{N}},
\; N \in \nnN,
\label{sofn}
\stboth
\\
C_{N}(x) & := & \sum_{n=0}^{\infty}  \frac{(-1)^{n} \cos(2n+1)x}{(2n+1)^{N}},
\; N \in \nnN.
\label{cofn}
\stboth
\ea
\no
The Hurwitz transform of $S_{N}$ and $C_{N}$ is expressed in terms of
Dirichlet's beta function
\ba
\beta(z) & := & \sum_{j=0}^{\infty} \frac{(-1)^{j}}{(2j+1)^{z}}
  =  4^{-z} \left[ \zeta(z, 1/4) - \zeta(z, 3/4) \right].
\stboth
\label{newbeta}
\ea
\no
Catalan's constant $G$ is $\beta(2)$. Properties of  $\beta(z)$ can be
found in \cite{atlas}, chapter 3.  \\

\no
\begin{Example}
Let $z \in \npR$. Then
\ba
\ione S_{N}(2 \pi q) \zeta(z,q) dq & = &
\frac{\Gamma(1-z)\cos(\pi z/2)}{(2 \pi)^{1-z}}
\beta(1-z+N),
\label{beta1}
\stboth
\\
\ione C_{N}(2 \pi q) \zeta(z,q) dq & = &
\frac{\Gamma(1-z)\sin(\pi z/2)}{(2 \pi)^{1-z}}
\beta(1-z+N).
\label{beta2}
\stboth
\ea
\end{Example}

\begin{proof}
The usual technique yields
\ba
\ione S_{N}(2 \pi q) \zeta(z,q) dq & = &
\frac{\Gamma(1-z)}{(2 \pi)^{1-z}} \cos(\pi z/2)
\sum_{n=0}^{\infty} \frac{(-1)^{n}}{(2n+1)^{1-z+N}}, \nn
\stboth
\ea
\no
which is (\ref{beta1}). The proof of (\ref{beta2}) is similar.  \\
\end{proof}

\medskip

The proof of the next two examples is similar to that of Example $10.1$.

\no
\begin{Example}
Let $m, N  \in \nnN$. Then
\ba
 & &  \\
\ione B_{m}(q) S_{N}(2 \pi q) dq & = &
\begin{cases}
0 & \text{ if } m \text{ is even,} \\
(-1)^{(m+1)/2} (2 \pi)^{-m} \, m! \, \beta(m+N) & \text{ if } m
\text{ is  odd,}
\end{cases}
\nn
\stboth
\ea
\no
and
\ba
 & & \\
\ione B_{m}(q) C_{N}(2 \pi q) dq & = &
\begin{cases}
(-1)^{m/2+1} (2 \pi)^{-m} \, m! \, \beta(m+N)  & \text{ if } m
\text{ is even,}  \\
0 & \text{ if } m \text{ is odd}.
\end{cases}
\nn
\stboth
\ea
\end{Example}

\medskip

\no
\begin{Example}
Let $m, N \in \mathbb{N}$. Then
\ba
\ione q^{m} S_{N}(2 \pi q) dq & = &
m!\sum_{k=0}^{\lfloor{\frac{m-1}{2} \rfloor}}
(-1)^{k+1} \frac{(2k+1)! \, \beta(2k+1+N)}{(m-2k)!(2 \pi)^{2k+1}}
\label{qsn}
\stboth
\ea
\no
and
\ba
\ione q^{m} C_{N}(2 \pi q) dq & = &
m! \sum_{k=0}^{\lfloor{\frac{m}{2} \rfloor}}
(-1)^{k+1}\frac{(2k)! \, \beta(2k+N)}{(m+1-2k)!(2 \pi)^{2k}}.
\nn
\stboth
\ea
\end{Example}

\medskip

\no
{\bf Note}. The values of $\beta$ at odd integers are given
by
\ba
\beta(2k+1) & = &
\frac{\abs{E_{2k}} }{2 (2k)!} \left( \frac{\pi}{2} \right)^{2k+1},
\label{betaodd}
\stboth
\ea
\no
where $E_{2k}$ are the Euler numbers defined by the generating
function
\ba
\frac{1}{\cos t} & = & \sum_{n=0}^{\infty} \frac{(-1)^{n}}{(2n)!} E_{2n}
t^{2n}. \label{eulergen}
\stboth
\ea
\no
We thus have
\begin{align}
\ione q^{m} S_{2n}(2 \pi q) \, dq& = \frac{m! \pi^{2n}}{2^{2n+3}}
\sum_{k=0}^{\lfloor\frac{m-1}{2} \rfloor}
\frac{(-1)^{k+1} (2k+1)! \, | E_{2n+2k} | }{2^{4k} \, (m-2k)! \, (2n+2k)! }
\stboth\\
\intertext{and}
\ione q^{m} C_{2n+1}(2 \pi q) \, dq& = \frac{m! \pi^{2n+1}}{2^{2n+2}}
\sum_{k=0}^{\lfloor\frac{m}{2} \rfloor}
\frac{(-1)^{k+1} (2k)! \, | E_{2n+2k} | }{2^{4k} \, (m-2k+1)! \, (2n+2k)! }.
\stoli
\end{align}

\medskip

\section{Eisenstein series}

In this section we compute the Hurwitz transform of functions related to
the Eisenstein series $G_{k}(\tau)$. \\

The Eisenstein series defined by\footnote{The sum is over $\mathbb{Z}^{2} -
(0,0)$.}
\ba
G_{k}(\tau) & := & \sum\nolimits_{m,n}^{'} \frac{1}{(m \tau + n)^{2k} },
\label{eisen}
\ea
\no
for $k \ge 2$ and $\tau \in \mathbb{C}$ with $\imagpart\tau > 0$,
are periodic functions with expansion (\cite{serre}, page 92):
\ba
G_{k}(\tau) & = & 2 \zeta(2k) + 2 \frac{(-1)^{k} (2 \pi)^{2k} }{(2k-1)!}
\sum_{n=1}^{\infty} \sigma_{2k-1}(n) e^{2 \pi i n \tau}, \label{expein}
\ea
\no
where
\ba
\sigma_{s}(n) & := & \sum_{d | n} d^{s}. \label{divisor}
\ea
\no
These series appear as coefficients in the cubic
\ba
y^{2} & = & 4x^{3} - 60G_{2}x - 140G_{3}
\ea
\no
that represents the torus
$\mathbb{C}/\mathbb{L}$, with $\mathbb{L}
= \mathbb{Z} \oplus \tau \mathbb{Z}$. See \cite{mcmoll} for details.  \\

Write $\tau = q + i t$, with $t >0$ and $0 \leq q \leq 1$. The
expansion (\ref{expein}) becomes
\begin{equation}
G_{k}(q + it)  =  2 \zeta(2k) + 2 \frac{(-1)^{k} (2 \pi)^{2k} }{(2k-1)!}
\sum_{n=1}^{\infty}\sigma_{2k-1}(n) e^{-2 \pi n t}\!
\left( \cos( 2 \pi k q) + i \sin(2 \pi k q) \right), \nn
\end{equation}
\no
so the Fourier coefficients of $G_{k}(q + i t)$ are
\ba
a_{n}  =
2 \frac{(-1)^{k} (2 \pi)^{2k} }{(2k-1)!}
\sigma_{2k-1}(n) e^{-2 \pi n t}   & \text{ and } & b_{n} = i a_{n}.
\label{general}
\ea
\no
We were unable to evaluate the corresponding Dirichlet series arising from
(\ref{general}). Instead we consider the functions
\ba
G_{k}^{(\alpha)}(q) & := &
\ift t^\alpha \left[ G_{k}(q+it) - 2 \zeta(2k) \right]
dt,
\label{modified}
\ea
\no
where $\alpha \in \mathbb{R}^{+}$.  We then have the following result. \\

\no
\begin{Example}
The Hurwitz transform of $G_{k}^{(\alpha)}(q)$ for $\alpha > z+2k-2$ is
\begin{multline}
\ione G_{k}^{(\alpha)}(q) \zeta(z,q) dq =
2\pi i\frac{{e^{ - i\pi z/2} }}{{\sin (\pi (\alpha  -
z)/2)}}\frac{{\Gamma (\alpha  + 1)\Gamma (1 - z)}}{{\Gamma
(2k)\Gamma (3 + \alpha  - z - 2k)}}\\
\times\zeta (2 + \alpha  - z)\zeta (3
+ \alpha  - z - 2k),
\label{inteins}
\end{multline}
where $z\in\npR$.
\no
\end{Example}
\begin{proof} The Fourier coefficients of $G_{k}^{(\alpha)}(q)$ are
\ba
a_{n}  =
2\Gamma(\alpha+1) \frac{(-1)^{k} (2 \pi)^{2k-\alpha-1} }{(2k-1)!}
\frac{\sigma_{2k-1}(n)}{n^{\alpha+1}}  & \text{ and } & b_{n} = i a_{n}.
\label{general1}
\ea
\no
The main theorem then yields
\ba
\ione G_{k}^{(\alpha)}(q) \zeta(z,q) dq & = &
\frac{2 \Gamma(\alpha+1)\Gamma(1-z) (-1)^{k} i
e^{-i\pi z/2}} {(2k-1)!(2 \pi)^{2-z-2k+\alpha}} \sum_{n=1}^{\infty}
\frac{\sigma_{2k-1}(n)}{n^{2-z+\alpha}}. \nn
\ea
\no
The last series is identified in \cite{apos}, page 231, as
\ba
\sum_{n=1}^{\infty} \frac{\sigma_{p}(n)}{n^{s}} & = & \zeta(s) \zeta(s -
p),\quad \realpart{s}>\max\{1,1+\realpart{p}\} \label{apostol}
\ea
\no
which, for $2-z+\alpha>2k$, implies (\ref{inteins}) after using
Riemann's relation (\ref{Riemann1}).\\
\end{proof}

\medskip

\section{Integrals over $[0, \tfrac{1}{2}]$. } \label{inthalf}

In this section we construct some examples of definite integrals over the
interval $[ 0 , \tfrac{1}{2}]$. Some of the examples described here are special
cases of the indefinite integral given at the end of Section \ref{examples}. \\

\medskip

\no
\begin{Example}
Let $ z \in \npR$. Then
\ba
\int_{0}^{1/2} \zeta(z,q) dq  &  = &
\frac{4 \Gamma(1-z) }{(2 \pi)^{2-z}} \cos\left( \frac{\pi z}{2} \right)
(1 - 2^{z-2}) \zeta(2 -z)
\label{half0}
\stboth
\\
 & = & \frac{2(2^{z-2}-1)}{1-z} \, \zeta(z-1).
\label{half01}
\stboth
\ea
\end{Example}

\begin{proof}
The Fourier expansion
\ba
\sum_{n=1}^{\infty} \frac{\sin 2 \pi(2 n+1) q}{2n+1} & = &
\begin{cases}
  \phantom{-}\tfrac{\pi}{4} &  \text{ if } \; 0 \leq q < \tfrac{1}{2}, \\
  -\tfrac{\pi}{4} &  \text{ if } \; \tfrac{1}{2} \leq q < 1
\end{cases}
\stboth
\ea
\no
yields, according to (\ref{riezeta1}),
\ba
\int_{0}^{1/2} \zeta(z,q) dq -
\int_{1/2}^{1} \zeta(z,q) dq  & = &
\frac{4}{\pi}
\frac{\Gamma(1-z) }{(2 \pi)^{1-z}} \cos \left( \frac{\pi z}{2} \right)
(1 - 2^{z-2}) \zeta(2 -z). \nn
\stboth
\ea
\no
But the vanishing of the integral of $\zeta(z,q)$ over $[0,1]$ can be written
as
\ba
\int_{0}^{1/2} \zeta(z,q) dq +
\int_{1/2}^{1} \zeta(z,q) dq  & = & 0, \nn
\stboth
\ea
\no
so (\ref{half0}) is proved.  The expression (\ref{half01})  follows
from (\ref{half0}) and Riemann's functional equation for the
$\zeta$-function. Alternatively, (\ref{half01}) can be derived directly
from the indefinite integral (\ref{intind}) and (\ref{riezeta1}).
\end{proof}

\medskip

\no
\begin{Example}
Let $m \in \nnN$. Then
\ba\label{beronehalf}
\int_{0}^{1/2} B_{m}(q) dq & = &
\frac{(-1)^{m+1} (2^{m+1}-1) B_{m+1}}{2^{m}(m+1)}.
\stboth
\ea
\no
\end{Example}
\begin{proof}
For $m\in\allN$ let $z= 1-m$ in (\ref{half01}) and
use
(\ref{berzeta1}). Formula (\ref{beronehalf}) can also be checked
to hold for the case $m=0$, using the value $B_1=-1/2$.
\end{proof}

\medskip

\begin{Thm}
The integrals
\ba
I(z,n) & := & \int_{0}^{1/2} q^{n} \zeta(z,q) dq, \quad n\in\allN
\label{momhalf}
\stboth
\ea
\no
satisfy the recursion relation
\ba
I(z,n) & = &
\frac{2^{z-1}-1}{2^{n}(1-z)} \zeta(z-1)
- \frac{n}{1-z} I(z-1,n-1).
\label{recul}
\stboth
\ea
\end{Thm}
\begin{proof}
Integrate by parts the identity
\ba
I(z,n) & = & \frac{1}{1-z} \int_{0}^{1/2} q^{n} \frac{\partial}{\partial q}
 \zeta(z-1,q) \, dq  \nn
\stboth
\ea
\no
and use
\ba
\zeta(z-1, \tfrac{1}{2}) & = & (2^{z-1} -1 ) \zeta(z-1) \label{zetahalf1}
\stboth
\ea
\no
to simplify the boundary terms. \\
\end{proof}

\no
Using the result (\ref{half01}) for $I(z,0)$, the recursion
relation yields the values
\begin{align}
I(z,1) & =  - \frac{(2^{z}-2) \zeta(z-1)}{4(z-1)} -
\frac{(2^{z}-8) \, \zeta(z-2)}{4(z-1)(z-2)}, \nn \\
\begin{split}
I(z,2) & =  - \frac{(2^{z}-2) \, \zeta(z-1)}{8(z-1)} -
\frac{(2^{z}-4) \, \zeta(z-2) }{4(z-1)(z-2)}
 - \frac{(2^{z}-16) \zeta(z-3)}{4(z-1)(z-2)(z-3)}. \nn
\end{split}
\end{align}
\medskip

A direct consequence of (\ref{recul}) is the following result.  \\

\begin{Thm}
\label{halfmomentsofzeta}
For $n \in \nnN$  and $z \in \npR$ let
\ba
I(z,n) & = & \int_{0}^{1/2} q^{n} \zeta(z,q) dq.  \nn
\stboth
\ea
\no
Then
\ba
I(z,n) =
n! \, \sum_{j=1}^{n+1} \frac{(-1)^{j} (1 - 2^{z-j}) \zeta(z-j)}
{(n+1-j)! 2^{n+1-j} (1-z)_{j}} +
(-1)^{n+1} n! \frac{\zeta(z-n-1)}{(1-z)_{n+1}}. \nn
\stboth
\ea
\end{Thm}

\no
\begin{Example}
Let $m,n \in \mathbb{N}$. Then
\begin{multline}
\int_{0}^{1/2} q^{n} B_{m}(q) dq =( - 1)^m \frac{{n!m!}}{{(n +
m + 1)!}}\\
\times\left[ {B_{n + m + 1}  + \sum\limits_{j = 1}^{n + 1}
{\binom{{n + m + 1}}{{n + 1 - j}}\left( {\frac{1}{{2^{n + 1 - j} }}
- \frac{1}{{2^{n + m} }}} \right)B_{m + j} } } \right].
\stboth
\end{multline}
\end{Example}
\begin{proof}
Let $z = 1-m$ in Theorem \ref{halfmomentsofzeta}.
\end{proof}

\medskip
\begin{Example}
We can use the result in Theorem \ref{halfmomentsofzeta} to give a
proof of Gosper's formula (\ref{gosper1}). From $\ln\Gamma(q+1)=\ln
q+\ln\Gamma(q)$ and (\ref{loggamma0}) we obtain
\ba
\int_{0}^{1/2} \ln \Gamma(q+1) dq = - \frac{1}{2} - \frac{{\ln
2}}{2} + \frac{{\lp }}{2} + \left. {\frac{d}{{dz}}}
\right|_{z = 0} \int_0^{1/2} {\zeta (z,q)dq}.
\nn
\ea
Now use
\ba
\left. {\frac{d}{{dz}}} \right|_{z = 0} \left[ {\frac{{2\left(
{2^{z - 2}  - 1} \right)}}{{(1 - z)}}\zeta (z - 1)} \right] =
\frac{1}{8} - \frac{{\ln 2}}{{24}} - \frac{3}{2}\zeta '( - 1) \stboth
\nn
\ea
to evaluate the last term above and produce
\ba
\int_{0}^{1/2} \ln \Gamma(q+1) dq =- \frac{3}{8} - \frac{{13\ln
2}}{{24}} + \frac{{\lp }}{2} - \frac{3}{2}\zeta '( - 1). \stboth
\ea
Finally, use (\ref{Riemann3}) to express $\zeta '( - 1)$ as
\ba\label{zetaprimeofone}
\zeta '( - 1) = \frac{{\zeta '(2)}}{{2\pi ^2 }} -
\frac{1}{{12}}(2\lp +\gamma  - 1)
\ea
to obtain (\ref{gosper1}). The derivation of (\ref{gosper2})
will appear in \cite{bem}.
\end{Example}

\medskip

\no
\section{A trigonometric example}

In this section we compute the Hurwitz transform of powers of sine and
cosine.  \\

\no
\begin{Example}
Let $z \in \npR$ and $n \in \allN$. Then
\begin{equation}
\label{sin2nzeta}
\ione \sin^{2n}(\pi q) \, \zeta(z,q) dq =
\frac{\Gamma(1-z)}{(2 \pi)^{1-z} 2^{2n-1} } \, \sin \left( \frac{\pi z}{2}
\right) \; \sum_{k=1}^{n} \frac{(-1)^{k}}{k^{1-z}} \binom{2n}{n-k}.
\stboth
\end{equation}
\end{Example}

\medskip

\no
\begin{Example}
Let $m, \, n \in \nnN$. Then
\begin{align}
\ione B_{2m+1}(q) \sin^{2n}(\pi q) dq &= 0,
\stboth\\
\intertext{and}
\ione B_{2m}(q) \sin^{2n}(\pi q) dq &=
\frac{(-1)^{m+1} (2m)!}{2^{2n-1} (2 \pi)^{2m} } \sum_{k=1}^{n}
\frac{(-1)^{k} }{k^{2m}} \binom{2n}{n-k}.
\stboth
\end{align}
\no
\end{Example}
\medskip
\no
The proof is a direct consequence of results (\ref{fou1}) and
(\ref{fou2}) for the Fourier coefficients of the Hurwitz zeta
function, once we expand $\sin^{2n}(\pi q)$ using a formula of Kogan
\cite{kogan}
\ba
\sin^{2n}x & = & \frac{1}{2^{2n}}  \binom{2n}{n} +
\frac{(-1)^{n}}{2^{2n-1}}
\sum_{k=0}^{n-1} (-1)^{k} \binom{2n}{k} \cos\left[2(n-k)x\right].
\stboth
\label{kogan1}
\ea

\medskip

Similar formulae exist for other powers of sine and cosine. Indeed,
\cite{kogan} shows
\ba
\sin^{2n+1}x & = & \frac{(-1)^{n}}{2^{2n}} \sum_{k=0}^{n} (-1)^{k}
\binom{2n+1}{k} \sin\left[(2n+1-2k)x \right], \nn
\stboth
\\
\cos^{2n}x & = & \frac{1}{2^{2n}} \binom{2n}{n} +
\frac{1}{2^{2n-1}}
\sum_{k=0}^{n-1} \binom{2n}{k} \cos\left[2(n-k)x\right], \nn
\stboth
\\
\cos^{2n+1}x & = &
\frac{1}{2^{2n}}
\sum_{k=0}^{n} \binom{2n+1}{k} \cos\left[(2n+1-2k)x\right], \nn
\stboth
\ea
\no
which yield
\begin{align}
\ione \sin^{2n+1}(2\pi q) \, \zeta(z,q) dq & =
\frac{\Gamma(1-z)}{(2 \pi)^{1-z} 2^{2n} } \, \cos \left( \frac{\pi z}{2}
\right) \; \sum_{k=0}^{n}  \frac{(-1)^k}{(2k+1)^{1-z}}\binom{2n+1}{n-k}, \nn
\stboth
\nn \\
\ione \cos^{2n}(\pi q) \, \zeta(z,q) dq & =
\frac{\Gamma(1-z)}{(2 \pi)^{1-z} 2^{2n-1} } \, \sin \left( \frac{\pi z}{2}
\right) \; \sum_{k=1}^{n} \frac{1}{k^{1-z}}\binom{2n}{n-k}, \nn
\stboth
\nn \\
\ione \cos^{2n+1}(2\pi q) \, \zeta(z,q) dq & =
\frac{\Gamma(1-z)}{(2 \pi)^{1-z} 2^{2n} } \, \sin \left( \frac{\pi z}{2}
\right) \; \sum_{k=0}^{n}  \frac{1}{(2k+1)^{1-z}}\binom{2n+1}{n-k}, \nn
\stboth
\nn
\intertext{and also}
\ione B_{2m+1}(q) \sin^{2n+1}(2\pi q) dq & =
\frac{(-1)^{m+1} (2m+1)!}{2^{2n} (2 \pi)^{2m+1} } \sum_{k=0}^{n}
\frac{(-1)^{k}}{(2k+1)^{2m+1}} \binom{2n+1}{n-k},
\stboth
\nn\\
\ione B_{2m}(q) \sin^{2n+1}(2 \pi q) dq & =  0, \nn \\
\ione B_{2m}(q) \cos^{2n}(\pi q) dq & =
\frac{(-1)^{m+1} (2m)!}{2^{2n-1} (2 \pi)^{2m} } \sum_{k=1}^{n}
\frac{1}{k^{2m}} \binom{2n}{n-k},
\stboth
\nn \\
\ione B_{2m+1}(q) \cos^{2n}(\pi q) dq & =  0, \nn \\
\nn\\
\ione B_{2m}(q) \cos^{2n+1}(2\pi q) dq & =
\frac{(-1)^{m+1} (2m)!}{2^{2n} (2 \pi)^{2m} } \sum_{k=0}^{n}
\frac{1}{(2k+1)^{2m}} \binom{2n+1}{n-k},
\stboth
\nn \\
\ione B_{2m+1}(q) \cos^{2n+1}(2 \pi q) dq & = 0. \nn
\end{align}

\begin{Example}
For $n\in\allN$
\begin{align}
\ione \sin ^{2n} (\pi q)\ln \Gamma (q)\,dq &= \frac{1}{{2^{2n + 1}
}}\sum\limits_{k = 1}^n {\frac{{( - 1)^k }}{k}\binom{{2n}}{{n - k}}
+ } \frac{1}{{2^{2n} }}\binom{{2n}}{n}\lp.
\stoli
\end{align}
\end{Example}
\begin{proof}
Simply use (\ref{loggamma0}), (\ref{sin2nzeta}) and Wallis' formula
\ba
\ione \sin ^{2n} (\pi q)\,dq = \frac{1}{{2^{2n} }}\binom{{2n}}{n}.
\stoli
\ea
\end{proof}

\section{The case $z$ positive}
\label{extzneg}

In this  section we extend some of the previous formulae to the
case $z \in \mathbb{R}^{+}$. Although the formulae of the previous
sections were derived under the assumption $ z \le 0$, so that
the Fourier expansion (\ref{fouzeta}) could be used, they can be
analytically extended to those positive values of $z$ where the
integral in question converges. This is so because the Hurwitz
transform (\ref{hurtrans1}) defines an analytic function of $z$ as
long as the defining integral converges.
For $z>0$ the only singularity of $\zeta(z,q)$ in the range
$0\le q\le 1$ lies actually at $q=0$, where it behaves as
$q^{-z}$. In fact,
\ba
\zeta(z,q) = \frac{1}{q^{z}} + \zeta(z,q+1),
\label{zetashift}
\ea
with $\zeta(z,q)$ finite for $q\ge 1$. The relation
(\ref{zetashift}) follows directly from the definition
(\ref{hurdef}) of the Hurwitz zeta function when $\realpart{z}>1$
and can be extended to the whole punctured complex $z$-plane,
$\mathbb{C}-\{1\}$, for $q>0$.\\

\begin{Example}
The formula (\ref{intber1}) derived in Example \ref{ex-berzeta}, namely
\ba
\ione B_{m}(q) \zeta(z,q) dq & = &
(-1)^{m+1} \frac{ m! \, \zeta(z-m) }{(1-z)_{m}},
\stboth
\nn
\ea
holds for real $z<1$ if $m$ equals one or an even integer, and for
$z<2$ otherwise.
\end{Example}
\begin{proof}
>From (\ref{polyber}) it is seen that near $q=0$ the Bernoulli
polynomials behave as
\ba
B_{m}(q)=B_{m}+ m B_{m-1}q+O(q^2),\quad m\ge 1.
\nn
\ea
Thus, the integrand $B_{m}(q) \zeta(z,q)$ behaves as $q^{-z}$ or
$q^{1-z}$, according if $B_m\ne 0$ or not. The result
now follows from the fact that the singularity $q^{-\alpha}$ is
integrable for $0<\alpha <1$.
\end{proof}

\begin{Example}
The formula (\ref{zetasq}) derived in Example \ref{ex-zetasq}, namely
\ba
\ione \zeta^{2}(z,q) dq & =
2 \Gamma^{2}(1-z) ( 2 \pi)^{2z-2} \zeta(2 - 2z)
\stboth
\nn
\ea
holds for real $z<1/2$.
\end{Example}
\begin{proof}
This follows directly from (\ref{zetashift}) and a reasoning
similar to the proof of the previous example.
\end{proof}

\medskip

In the rest of this section use integration by parts to derive
formulas for integrals containing the function
\ba
\zeta_{*}(z,q)  := \zeta(z,q+1)=
\zeta(z,q) - \frac{1}{q^{z}}, \label{zetastar}
\ea
which is finite in the closed interval $[0,1]$ for arbitrary $z
\neq 1$.

\begin{Thm}
Let $f$ be $n$-times differentiable and $z\in\allR
-\{1,2,\ldots,n+1\}$. Then
\begin{equation} \label{zetaparts}
\begin{split}
\ione f(q)  \zeta_{*}(z,q) dq &  =    - \sum_{k=0}^{n-1} (-1)^{k}
\frac{f^{(k)}(1)}{(1-z)_{k+1}} \\ & + \sum_{k=0}^{n-1} (-1)^{k}
\frac{f^{(k)}(1) - f^{(k)}(0) }{(1-z)_{k+1}} \, \zeta(z-k-1)  \\
 &  +  \frac{(-1)^{n}}{(1-z)_{n}} \ione f^{(n)}(q) \zeta_{*}(z-n,q) dq. \\
\end{split}
\stboth
\end{equation}
\end{Thm}
\begin{proof}
Integrate by parts to produce
\ba
\ione f(q) \zeta_{*}(z,q) dq & = & \frac{f(1)-f(0)}{1-z}
\zeta(z-1)  - \frac{f(1)}{1-z} - \label{zetaone}  \\
 & - & \frac{1}{1-z} \ione f'(q) \zeta_{*}(z-1,q) dq,
\stboth
\nn
\ea
\no
which establishes the result for $n=1$.
Repeated integration by parts yields (\ref{zetaparts}).
\end{proof}

We now apply this theorem to a few well chosen functions $f$.

\medskip
\no
\begin{Example}
Take $f(q)=1$ and $n=1$ in
(\ref{zetaparts}). We get
\ba
\int_0^1 {\zeta _* (z,q)\,dq}  =
\frac{1}{{z - 1}}  \label{zetastaralone}
\ea
\no
for $z \neq 1$.
\end{Example}
\medskip
\no
\begin{Thm}
Let $r\in \nnN$ and $z\in\allR -\{1,2,\ldots,2r+2\}$. Then
\begin{multline}
\int_0^1 {\zeta _* (z - 2r - 1,q)\,\zeta _* (z,q)\,dq}  =
\frac{1}{{2\,(1 + r - z)}}\\ - \frac{1}{2}\sum\limits_{k = 0}^{2r}
{\frac{{(z - 2r - 1)_k }}{{(1 - z)_{k + 1} }}\left[ {\zeta (z - k
- 1) + \zeta (z - 2r - 1 + k)} \right]}. \label{zetastar0}
\end{multline}
\end{Thm}
\begin{proof}
Take $f(q)=\zeta_*(z',q)$ in (\ref{zetaparts}). We obtain
\[
f^{(k)} (q) = ( - 1)^k (z')_k \zeta _* (z' + k,q),
\]
so that
\begin{align}
\nn f^{(k)} (1)& = ( - 1)^k (z')_k \left[ {\zeta (z' + k) - 1}
\right] \\ \intertext{and} \nn f^{(k)} (0)& = ( - 1)^k (z')_k
\zeta (z' + k).
\end{align}
Hence
\begin{multline}
\int_0^1 {\zeta _* (z',q)\,\zeta _* (z,q)\,dq}  = \frac{{(z')_n
}}{{(1 - z)_n }} \int_0^1 {\zeta _* (z'+n,q)\,\zeta _*
(z-n,q)\,dq}\\ - \sum\limits_{k = 0}^{n-1} {\frac{{(z')_k }}{{(1 -
z)_{k + 1} }}\left[ {\zeta (z - k - 1) + \zeta (z' + k) -1}
\right]}. \label{zetastar1}
\end{multline}
Now choose $z'=z-n$ with  $n=2r+1$. Then
\[
(z')_n  = (z - (2r + 1))_{2r + 1}  = ( - 1)^{2r + 1} (1 - z)_{2r +
1}  =  - (1 - z)_{2r + 1},
\]
in virtue of the identity $(z-j)_j=(-1)^j(1-z)_j$. Moreover,
$z'+n=z$ and $z-n=z'$. Thus the integral on the right hand side of
(\ref{zetastar1}) is just the negative of the initial integral. The
result (\ref{zetastar0}) follows, since \ba
\sum_{k=0}^{2r}{\frac{{(z - 2r - 1)_k }}{{(1 - z)_{k + 1} }} =
\frac{1}{{1 + r - z}}}. \label{pochhammersum1} \ea
\end{proof}

\medskip
\no
\begin{Example}
The case $r=0$ in (\ref{zetastar0}) yields
\ba \int_0^1 {\zeta _* (z - 1,q)\,\zeta _* (z,q)\,dq}  =
\frac{{2\zeta (z - 1) - 1}}{{2\,(z - 1)}}. \label{zetastar0withr0}
\ea
\no
This result can be obtained alternatively by noting that
the integrand has a simple antiderivative, namely,
\ba
\int {\zeta
_* (z - 1,q)\,\zeta _* (z,q)\,dq =  - \frac{1}{{2\,(z - 1)}}\zeta
_*^2 (z - 1,q)}. \label{antideriv2}
\ea
\no
The expression (\ref{zetastar0withr0})
now follows from the evaluation
\[
\left. {\zeta _*^2 (z - 1,q)} \right|_0^1  = 1 - 2\zeta (z - 1).
\]
\end{Example}
\medskip
\no
\begin{Example}
The cases $\{z=5,r=1\}$ and
$\{z=5/2,r=2\}$ yield, respectively,
\ba
\int_0^1 {\zeta _*
(2,q)\,\zeta _* (5,q)\,dq}  =  - \frac{1}{6} + \frac{{\pi ^2
}}{{24}} + \frac{{\pi ^4 }}{{360}} - \frac{\zeta(3)}{6}
\label{intz2and5}
\stboth
\ea
and
\begin{multline}
\int_0^1 {\zeta _* ( - \tfrac{5}{2},q)\,\zeta _* (\tfrac{5}{2},q)\,dq}  = 1 +
\tfrac{2}{3}\zeta (- \tfrac{5}{2}) + \tfrac{{10}}{3}\zeta ( - \tfrac{3}{2})\\ -
10\zeta ( - \tfrac{1}{2}) + \tfrac{10}{3}\zeta (\tfrac{1}{2}) +
\tfrac{2}{3}\zeta( \tfrac{3}{2}).
\stboth
\end{multline}
\end{Example}

\medskip
\no
\begin{Thm}
Let $n\in\nnN$ and $z\in\allR -\{1,2,\ldots,n+1\}$. Then
\begin{equation}
\int_0^1 {q^n\,\zeta _* (z,q)\,dq}  = - \frac{1}{{1 - z + n}} +
n!\sum\limits_{k = 0}^{n - 1} {( - 1)^k \frac{{\zeta (z - k -
1)}}{{(n - k)!(1 - z)_{k + 1} }}}. \label{zetastarqn}
\end{equation}
\end{Thm}
\begin{proof}
Take $f(q)=q^n$ in (\ref{zetaparts}). Then
\[
f^{(k)} (q) = \frac{{n!}}{{(n - k)!}}q^{n - k} \;\text{for} \;k
\le n,\quad\text{and}\quad f^{(k)} (q) = 0\;\text{for} \;k > n.
\]
Thus
\begin{alignat}{3}
\nn f^{(k)} (1) &= \frac{{n!}}{{(n - k)!}}&\;\text{for} \;k <
n,&\quad\text{and}\quad f^{(n)} (1) &= 1,\\ \nn f^{(k)} (0) &= 0
&\;\text{for} \;k < n,&\quad\text{and}\quad f^{(n)} (0)& = 1.
\end{alignat}
Using (\ref{zetastaralone}) we have
\begin{equation}
\int_0^1 {q^n\,\zeta _* (z,q)\,dq}  = n!\sum\limits_{k = 0}^{n -
1} {( - 1)^k \frac{{\zeta (z - k - 1) - 1}}{{(n - k)!(1 - z)_{k +
1} }}}  - \frac{{( - 1)^n n!}}{{(1 - z)_{n + 1} }}. \nn
\end{equation}
The result now follows after using the identity
\begin{equation}
n!\sum\limits_{k = 0}^n {\frac{{( - 1)^k }}{{(n - k)!(1 - z)_{k +
1} }}}  = \frac{1}{{1 - z + n}}. \label{pochhammersum2}
\end{equation}
\end{proof}

\medskip
\no
For $n-z+1>0$ we have $\int_0^1 {q^{n - z} dq} = (1-z+n)^{-1}$, so
(\ref{zetastarqn}) yields
\[
\int_0^1 {q^n \,\zeta (z,q)\,dq}  = n!\sum\limits_{k = 0}^{n - 1}
{( - 1)^k \frac{{\zeta (z - k - 1)}}{{(n - k)!(1 - z)_{k + 1} }}},
\]
which is exactly of the same form as the result (\ref{mom})
derived for the case $z\le 0$. We can therefore combine both
results as the following.

\begin{Thm}
\label{thm-mom2}
Let $n\in\nnN$ and $z\in\allR$ such that $n-z+1>0$. Then, \ba
\int_0^1 {q^n \,\zeta (z,q)\,dq}  = n!\sum\limits_{k = 0}^{n - 1}
{( - 1)^k \frac{{\zeta (z - k - 1)}}{{(n - k)!(1 - z)_{k + 1} }}}.
\label{mom2} \ea
\end{Thm}

\section{Polygamma functions}

In this last section we evaluate the moments of the
polygamma functions, defined as
\ba \psi^{(n)}(z) & := & \frac{d^{n}}{dz^{n}} \psi(z),\quad
n\in\nnN,
\label{polyga}
\ea
\no
where $\psi^{(0)}(z)=\psi(z)$ is the digamma function defined in
(\ref{digamma}).\\

The polygamma functions can be expressed in terms of the
Hurwitz zeta function as (see \cite{atlas}, chapter $44$)
\begin{align}
\psi ^{(m)} (q) &= ( - 1)^{m + 1} m!\zeta (m + 1,q)\qquad
m=1,2,\ldots
\label{poly-hurwitz} \\
\intertext{and}
\psi(q)&=\lim_{z\to1}\left[ \frac{1}{z-1}-\zeta(z,q)\right],
\label{psi-hurwitz}
\end{align}
so the integrals given in Section \ref{extzneg} reduce to
integrals involving $\psi^{(n)}(z)$ when the variable $z$ is
made to approach a positive integer.
\no
\begin{Thm}
Let $n,m\in\allN$ with $n>m$. Then
\begin{equation}\label{mompolys}
\begin{split}
\int_0^1 q^n \,\psi ^{(m)} (q)\,dq  = (-1)^m &\frac{n!}{(n -
m)!}\Bigg[ \frac{\gamma }{n - m + 1}
+ (n - m)!\sum\limits_{k = 0}^{m - 2} \frac{\Gamma (m - k)\zeta
(m - k)}{(n - k)!}\\
&+ \sum\limits_{k = 0}^{n - m - 1} (-1)^k \binom{n - m}{k}\left[
H_k \zeta ( - k) + \zeta '( - k) \right] \Bigg].\\
\end{split}
\end{equation}
\end{Thm}
\begin{proof}
We compute the limit as $z\to m+1$ in Theorem \ref{thm-mom2}. Substitute
$z=m+1-\varepsilon$ in (\ref{mom2}) and let $\varepsilon\to
0$. We encounter two types of singularities as $\varepsilon\to
0$: one corresponding to the pole of $\zeta(s)$ at $s=1$, for
$k=m-1$, and the other corresponding to the vanishing of
the Pochhammer symbol $(-m)_{k+1}$, for $k=m,m+1,\ldots,n-1$. To
derive (\ref{mompolys}) consider the Laurent expansion of
(\ref{mom2}) about $\varepsilon=0$ up to order $\varepsilon^0$.
The following expansions are employed:
\begin{align}
\zeta (1 - \varepsilon ) &= \frac{1}{\varepsilon } + \gamma  + O(\varepsilon
),\\
\zeta ( - r - \varepsilon ) &= \zeta ( - r) - \varepsilon \zeta '(
- r) + O(\varepsilon ), \quad r=0,1,2,\ldots, \nn \\
\Gamma (m + 1)\frac{{\Gamma (1 - \varepsilon )}}{{\Gamma (m + 1 -
\varepsilon )}} &= 1 + H_m \varepsilon  + O(\varepsilon ), \nn \\
\Gamma (m + 1)\frac{{\Gamma ( - r - \varepsilon )}}{{\Gamma (m + 1
- \varepsilon )}} &= \frac{{( - 1)^{r + 1} }}{{r!}}\left[
{\frac{1}{\varepsilon } + H_m  - H_r } \right] + O(\varepsilon ), \quad
r=0,1,2,\ldots. \nn
\end{align}
A direct calculation yields
\begin{align}
{( - 1)^k \frac{{\zeta (z - k - 1)}}{{(n - k)!(1 - z)_{k +
1} }}} \Bigg|
_{\substack{z = m + 1 - \varepsilon\\k = m - 1\hfill}}
&= \frac{1}{{m!(n - m-1)!}}\left[ {\frac{1}{\varepsilon }
+ H_m  - \gamma} \right], \nn \\
\intertext{and, for $r=0,1,\ldots,n-m-1$,}
{( - 1)^k \frac{{\zeta (z - k - 1)}}{{(n - k)!(1 - z)_{k +
1} }}} \Bigg|_{\substack{z = m + 1 - \varepsilon\\k = m
+r\hfill}}
&= \frac{{( - 1)^r }}{{m!r!(n - m + r)!}} \nn \\
&\times\left[
{\frac{{\zeta ( - r)}}{\varepsilon } + (H_m  - H_r )\,\zeta ( - r)
- \zeta '( - r)} \right].
\nn
\end{align}
The coefficient of the singular term $1/\varepsilon$ is
\ba
\frac{1}{{m!(n - m)!}}\left[ {\frac{1}{{n - m + 1}} +
\sum\limits_{r = 0}^{n - m - 1} {( - 1)^r \binom{{n - m}}{r}\zeta (
- r)} } \right],
\nn
\ea
which vanishes in view of the identity
\ba\label{ident1}
\sum\limits_{r = 0}^{j} {( - 1)^r \binom{{j+1}}{r}\zeta (
- r)}  =  - \frac{1}{{j+2}}.
\ea
The rest of the terms can
be collected to yield (\ref{mompolys}), after multiplying by the
overall factor $(-1)^{m+1} m!$ in
(\ref{poly-hurwitz}).
\end{proof}
\no
Along similar lines, we can use relation (\ref{psi-hurwitz}) to
prove the following result.

\begin{Thm}
For $n\in\allN$,
\begin{equation}\label{momdigamma}
\int_0^1 q^n \,\psi (q)\,dq  = \zeta'(0)   +
\sum\limits_{k = 1}^{n - 1} (-1)^k \binom{n}{k}\left[
H_k \zeta ( - k) + \zeta'(-k) \right].
\end{equation}
\end{Thm}
\begin{proof}
Theorem \ref{thm-mom2}  and (\ref{psi-hurwitz}) yield
\ba
\ione q^{n} \psi(q) dq & = & \lim_{z \to 1} \ione q^{n}
\left[ \frac{1}{z-1} - \zeta(z,q) \right] \, dq \nn \\
& = & \lim_{z \to 1} \frac{1}{z-1}
\left[ \frac{1}{n+1} + \zeta(z-1) + n! \sum_{k=1}^{n-1}
\frac{(-1)^{k} \zeta(z-k-1)}{(n-k)! (2-z)_{k}} \right]  \nn
\ea
\no
and (\ref{momdigamma}) follows by l'Hopital's rule.
\end{proof}

\no
\section{Conclusions}

We have evaluated a series of definite integrals whose integrands
involve the Hurwitz zeta function. \\

Most of the formulae involving elementary functions and
$\ln\Gamma(q)$ can be considered to be special cases of Theorem
\ref{thm-zetazeta}, in view of the relations (\ref{bernoulli}),
(\ref{loggamma0}) and (\ref{reflec}), the latter written in the
form\\
\ba
\ln \Gamma (q) + \ln \Gamma (1 - q) = \ln \pi  - \ln \sin \pi q,
\label{logsum}
\ea
\\
which respectively relate the Bernoulli polynomials, the logarithm of the gamma
function and thus also the function $\ln \sin \pi q$ to $\zeta(z,q)$ in a linear
way.\\

\no
{\bf Acknowledgments}. The authors would like to thanks G. Boros for many
suggestions. The first author would like to thank the Department of
Mathematics at Tulane University for its hospitality during a
short visit, where this work was begun, and the support of CONICYT
(Chile) under grant 1980149.

\end{document}

--=====================_115625940==_
Content-Type: text/plain; charset="us-ascii"; format=flowed

Olivier Espinosa (espinosa@fis.utfsm.cl)
--=====================_115625940==_--